\documentclass[10pt,twocolumn,twoside]{IEEEtran} 

\usepackage[]{footmisc}

\usepackage{mathrsfs}
\usepackage{amsthm}

\newtheorem{lemma}{Lemma}
\newtheorem{remark}{Remark}

\usepackage{hyperref} 

\usepackage{algorithm,algorithmic}

\usepackage{amsfonts}
\usepackage{xcolor}
\usepackage{caption}
\usepackage{graphicx}

\usepackage[utf8]{inputenc}

\usepackage{nomencl} 
\makenomenclature

\usepackage{etoolbox}
\renewcommand{\nomgroup}[1]{%
  \item[\bfseries
  \ifstrequal{#1}{A}{Given Parameters}{%
  \ifstrequal{#1}{B}{Decision variables}{%
  \ifstrequal{#1}{C}{Other Symbols}{}}}%
]}
\usepackage{cite}

\ifCLASSINFOpdf
\else
\fi
\usepackage{amsmath}
\usepackage{algorithmic}

\usepackage{array}

\usepackage{subfigure}

\usepackage{float}

\usepackage{stfloats}

\begin{document}


\title{Joint Routing and Charging Problem of Multiple Electric Vehicles: A Fast Optimization Algorithm }
%
%
%

\author{Canqi~Yao\thanks{This work was supported in part by the National Key Research and
Development Program of China under Grant 2019YFB1705401; in part by
the Natural Science Foundation of China under Grant 61873118 and Grant
61903179; in part by the Science, Technology and Innovation Commission
of Shenzhen Municipality under Grant ZDSYS20200811143601004 and Grant
RCBS20200714114918137\\ \indent\indent C. Yao is with the School of Mechatronics Engineering, Harbin Institute of Technology, Harbin, 150000, China, is with the Shenzhen Key Laboratory of Biomimetic Robotics and Intelligent Systems, Department of Mechanical and Energy Engineering, and the Guangdong Provincial Key Laboratory of Human-Augmentation and Rehabilitation Robotics in Universities, Southern University of Science and Technology, Shenzhen 518055, China  (e-mail: vulcanyao@gmail.com).}, Shibo Chen\thanks{S. Chen and Z. Yang are with the Shenzhen Key Laboratory of Biomimetic Robotics and Intelligent Systems, Department of Mechanical and Energy Engineering, and the Guangdong Provincial Key Laboratory of Human-Augmentation and Rehabilitation Robotics in Universities, Southern University of Science and Technology, Shenzhen 518055, China (e-mails: chensb@sustech.edu.cn, yangzy3@sustech.edu.cn).}, \textit{Member, IEEE,} and Zaiyue Yang, \textit{Member, IEEE}}

\markboth{Accepted by IEEE Transactions on Intelligent Transportation Systems, doi: 10.1109/TITS.2021.3076601}%
{Yao \MakeLowercase{\textit{et al.}} Joint Routing and Charging Problem of Multiple Electric Vehicles: A Fast Optimization Algorithm}

\maketitle

\begin{abstract}
Logistics has gained great attentions with the prosperous development of e-commerce, which is often seen as the classic optimal vehicle routing problem. Meanwhile, electric vehicle (EV) has been widely used in logistic fleet to curb the emission of green house gases in recent years. Thus, solving the optimization problem of joint routing and charging of multiple EVs is in a urgent need, whose objective function includes charging time, charging cost, EVs travel time, usage fees of EV and revenue from serving customers. This joint problem is formulated as a mixed integer programming (MIP) problem, which, however, is NP-hard due to integer restrictions and bilinear terms from the coupling between routing and charging decisions. The main contribution of this paper lies at proposing an efficient two-stage algorithm that can decompose the original MIP problem into two linear programming (LP) problems, by exploiting the exactness of LP relaxation and eliminating the coupled term. This algorithm can achieve a near-optimal solution in polynomial time. In addition, another variant algorithm is proposed based on the two-stage one, to further improve the quality of solution.   \textcolor{black}{Compared with the state-of-the-art algorithm, extensive simulations are implemented to validate the effectiveness of the proposed algorithm. } \end{abstract}


\begin{IEEEkeywords}
Electric vehicle, routing problem, mixed integer programming, linear programming
\end{IEEEkeywords}

%
\IEEEpeerreviewmaketitle


\makenomenclature
\nomenclature[A]{$g_i$}{Charging time per unit power(1 kWh) at node $i,i\in\mathcal{V}'$ }
\nomenclature[A]{$p_i$}{Charging price at node $i,i\in\mathcal{V}'$}
\nomenclature[A]{$d_{ij}$}{Travel distance from node $i$ to $j$ where $i,j\in\mathcal{V}'$}
\nomenclature[A]{$T_{ij}$}{Travel time from node $i$ to $j$ where $i,j\in\mathcal{V}'$}
\nomenclature[A]{$e_{ij}$}{Energy consumption traveling from node $i$ to node $j$, $e_{ij}=\epsilon d_{ij}$, where $\epsilon$ is the energy consumption rate}
\nomenclature[A]{$\mathscr{M}_{i}$}{ Revenue from serving requests when $i\in\mathcal{R}$}
\nomenclature[A]{$c_v$}{EVs' usage cost }
\nomenclature[A]{$t_i$}{Pickup time of transportation requests at node $i,i\in\mathcal{R}$ given by customers}
\nomenclature[A]{$\mathcal{K}$}{The set of EVs}
\nomenclature[A]{$\mathcal{R}$}{The set of transportation request
 nodes}

\nomenclature[B]{$x_{ij}^k$}{Binary decision variable indicating whether electric vehicle $k\in\mathcal{K}$ traverses road $(i,j)\in\mathcal{E}'$ or not}
\nomenclature[B]{$r_i^k$}{Continuous decision variable indicating charged energy at node $i,i\in\mathcal{V}'$}

\nomenclature[C]{$E^k_i$}{Battery energy level of $k,k\in\mathcal{K}$ at node $i,i\in\mathcal{V}'$}
\nomenclature[C]{$G(\mathcal{V},\mathcal{E})$}{Directed and connected graph, $\mathcal{V}=\left\{v_1,\ldots,v_n \right\}$ consisting of the start depot, end depot and road intersection represents the nodes within graph $\mathcal{G}$ and edges between nodes are denoted by $\mathcal{E}$.}
\nomenclature[C]{$G'(\mathcal{V}',\mathcal{E}')$}{Directed graph, $V'=\left\{v_1,v_n \right\}\cup\mathcal{R}$ includes the start depot, end depot and transportation request nodes.}
\nomenclature[C]{$\Phi^k_j$}{Indicating whether the actual pickup time of node $j$ satisfies predefined pickup time or not.}
\nomenclature[C]{$\mathcal{R}^q_v$}{The  violated pickup time transportation requests set at $q^{th}$ iteration} 



\printnomenclature

\section{Introduction}
%
%
%
%

\IEEEPARstart{L}{ogistics}  and transportation (L\&T) activities represent a key sector in worldwide economy, and \textcolor{black}{ are significant contributors} to social and economic progress in modern society. Especially with the development of e-commerce, there were 87 billion parcels shipped worldwide in 2018 in total, which constitutes an 104 percent increase compared to 2014\cite{parcel}. As a result, the number of vehicles has a significant increase in response to the rise of L\&T activities\cite{juan2016electric}.

\par Conventional vehicles, e.g. internal combustion engine vehicles are greenhouse gases producers with higher maintenance costs and operational cost than Electric Vehicles (EVs). To reduce its carbon emission, Amazon places an order for 100,000 EVs in a push to make the company’s fleet run entirely on renewable energy\cite{amazon}. It is foreseen that a large amount of transportation requests will be completed by EVs in the near future, which complicates the decision process.
\par The basic transportation request delivery problem is a kind of vehicle routing problem (VRP), which can be formulated as the integer programming (IP)\cite{toth2002vehicle}. However, with the introduction of EVs, both route selection variables (binary) and electricity charging amount variables (continuous) are coupled together and need to be optimized jointly, which makes it a mixed integer programming (MIP) that is NP-hard and not easy to solve\cite{tang2018distributed}. 

\par Recently, significant research efforts have been put forward to solve the joint routing and charging problem, and the methods can be classified into three categories: a) commercial optimization solvers; b) heuristics method; c) problem-tailored algorithm based on structure information. Some researchers formulate routing and charging problem of EVs as an MIP problem solved by commercial optimization solvers directly\cite{chen2018optimal,Cerna7862287,james}, which requires remarkable computation resources and is unaffordable for most of applications. To speed up the computation time, some heuristics techniques are employed, such as evolutionary algorithm\cite{yang2015electricGenetic,lam2016autonomous,abdulaal2017solving} or the neighborhood search algorithm\cite{schneider2014electric,keskin2016partial,sassi2015multi}. A modified custom genetic algorithm incorporated with embedded Markov decision process and trust region optimization methods is developed to solve multivariant EV Routing Problem\cite{abdulaal2017solving}. However, there is no guarantee on the solution optimality of heuristic methods.

\par The structure information of joint routing and charging problem of EVs, which can be used to design tailored algorithm to reduce the time complexity\cite{Pourazarm,tang2018distributed,james2018autonomous,pourazarm2016optimal}, is not fully exploited.  Recently, the dual decomposition method is used in \cite{james2018autonomous} to  solve the joint routing and charging problem formulated in the form of quadratic-constrained mixed integer linear program with subproblem also modelled as MIP problem, which still requires remarkable computation resources. Besides, \cite{tang2018distributed} decomposes routing and charging scheduling optimization for multiple EVs into two subproblems (path selection problem and charging scheduling problem) and solves them sequentially with Lagrangian relaxation method in distributed manner. However, the charging time is not considered in \cite{tang2018distributed}, which is the critical feature of EVs. In \cite{Pourazarm}, the single EV routing and charging problem is modelled as the mixed integer nonlinear program, then reformulated as two linear programming (LP) problems. However, as in practice a fleet of EVs will be sent out for delivery, efficient algorithms for joint routing and charging problem of multiple EVs considering charging time constraint are still in an urgent need.

\par To this end, in this paper we propose a novel method which can provide a near optimal solution within polynomial time by fully exploiting the structure information of routing and charging problem of multiple EVs. In particular, we consider a fleet of EVs serving  transportation requests in a network with charging facilities, and the goal is to determine an optimal routing and charging plan with the minimization of charging time, charging cost, EVs travel time, usage fees of EV and revenue from serving customers.

\par  The main contributions of this paper are summarized as follows:
\begin{itemize}
    \item We formulate the joint  routing and charging problem of multiple EVs considering time constraints for each customer as an MIP problem.
    \item To relieve computational burden of MIP problem, we propose a computational efficient two-stage algorithm, which decomposes the original MIP problem into two LP problems, by exploiting the exactness of LP relaxation and eliminating the coupled term. This algorithm can achieve a near-optimal solution in polynomial time.
    \item To further improve the quality of solution, another variant algorithm is proposed based on the two-stage one.
\end{itemize}
\par The paper is organized as follows. We elaborate the system components of EVs routing and charging problem in Section II. Then an MIP problem concerning joint routing and charging of multiple EVs is formulated in Section III. In Section IV we propose a computational efficient two-stage method to solve EVs routing and charging problem. The iterative optimization variant method is  described in Section V.  Extensive simulations are implemented in Section VI to validate the effectiveness of our algorithm. Finally, we conclude this paper in Section VII.


\section{System description}

In this section, we firstly introduce the model of joint routing and charging problem of multiple EVs and the operations of the system. Then the preprocessing step is described in Subsection B to reduce the number of nodes. In the following sections, nodes and transportation requests are used interchangeably.
\subsection{System Components}
\subsubsection{Transportation Network}
A fleet of EVs traverses in the transportation network for finishing the logistic requests while meeting the pickup time specified by customers. The transportation network is modeled as a directed graph $G(\mathcal{V},\mathcal{E})$, where $\mathcal{V}=\{v_1,\ldots,v_n\}$ represents the set of road intersections and $\mathcal{E}$ denotes the set of road segments with $(i,j)\in\mathcal{E}$ denoting a road link from vertex $i$ to $j$. The start depot and end depot are denoted as $v_1$ and $v_n$ respectively. We denote $d_{ij}$ and $T_{ij}$ as the distance and travel time of road segments between road intersections $i,j\in \mathcal{V}$ respectively.

\subsubsection{Transportation Requests}
A subset of vertex in graph $G$ represents transportation requests and is denoted as $\mathcal{R}$,$\mathcal{R}\subset \mathcal{V}$. In this paper, we use a tuple $\left\{t_i, \mathscr{M}_{i} \right\},i\in \mathcal{R}$ to describe the characteristics of transportation  requests. Here $\mathscr{M}_{i}$ is the revenue from serving requests when $i\subset \mathcal{R}$. At the beginning of the scheduling phase, customers will broadcast their desired pickup time $t_i,i\in\mathcal{R}$. In this paper, we assume that every transportation request vertex $i\in \mathcal{R}$ has charging facilities with charging duration per unit energy and price denoted by $g_i,p_i,i\in \mathcal{R}$ respectively. Note that if node $i$ does not have a charging facility, charging duration per unit energy $g_i$ is simply set to $\infty$.

\subsubsection{Electric Vehicles and Charging Process}
The set of all EVs is denoted as $\mathcal{K}$. Each EV $k\in \mathcal{K}$ can be described by a tuple $\left\{E_k^{0},E_{k,max},c_{v}\right\}$. $E_k^{0}$ represents the initial energy of $k$ starting from the depot, while the capacity of EV $k$ is denoted as $E_{k,max}$. $c_{v}$ is the EVs' usage fee. All EVs station at the starting depot and return to the end depot. In the joint routing and charging problem of EVs, each EV is assigned with a set of transportation requests to serve. With limited battery size, the owner of EV has to replenish EV with charging facility. The charging amount of EV $k$ at node $i$  is denoted by $r^k_i, i\in \mathcal{R}, k\in\mathcal{K}$. The maximum charging amount provided by charging facility is denoted by $r_{max}$.

\subsection{Preprocessing}

\par To reduce the number of nodes of considered graph, we adopt the preprocessing idea from \cite{lam2016autonomous} transforming $G(\mathcal{V},\mathcal{E})$ to a much smaller graph $G'(\mathcal{V}',\mathcal{E}')$. $\mathcal{V}'$ includes the start depot, transportation request nodes and end depot, i.e. $\mathcal{V}'=\{v_1,v_n\}\cup \mathcal{R}$ and $\mathcal{R}\subseteq \mathcal{V}\setminus\{v_1,v_n\}$. 
\par
\begin{figure*}[ht]
\centering
\includegraphics[width=.85\linewidth]{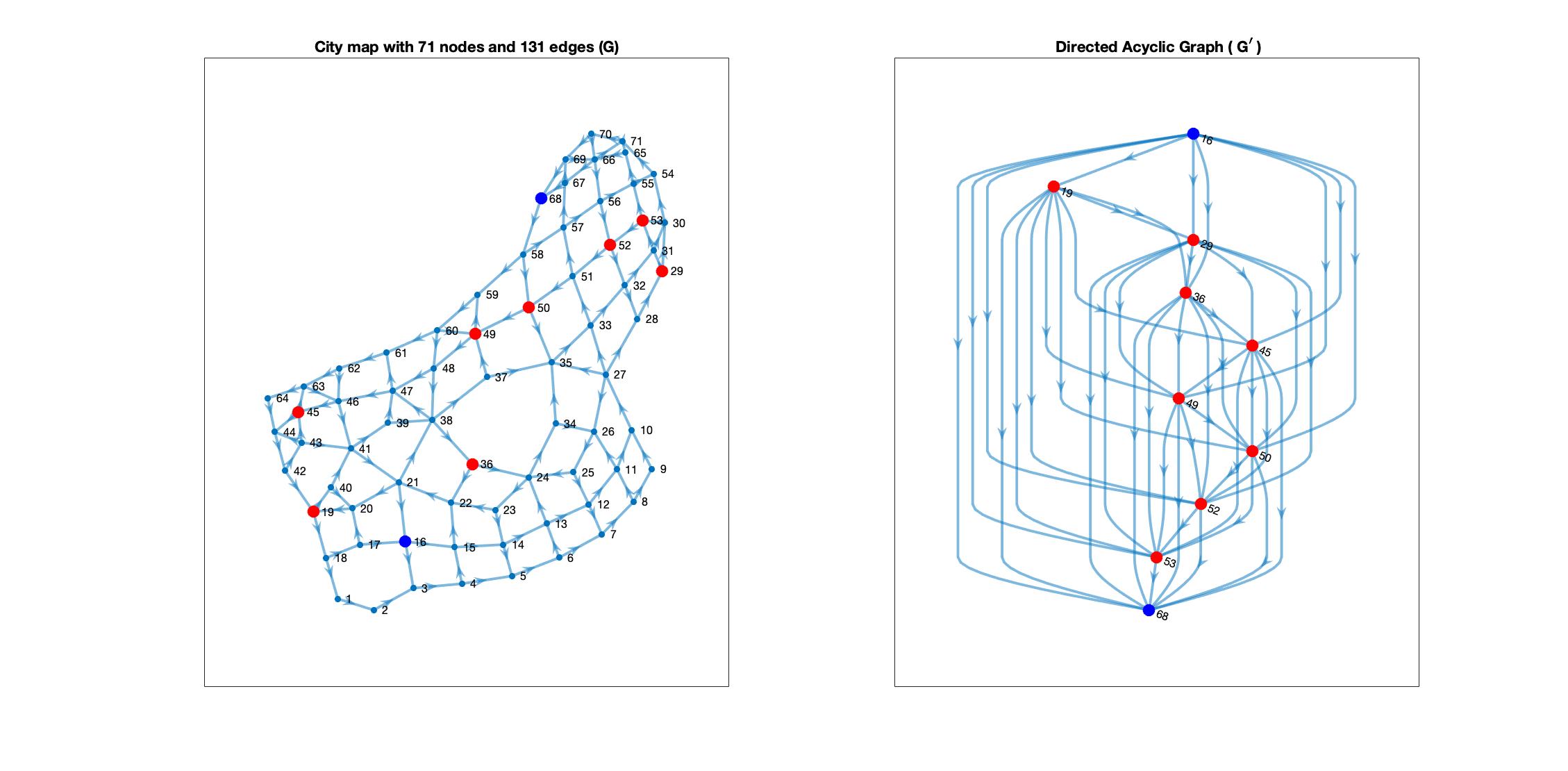}
  \caption{An example of real road network and $G^{'}$}
\label{roadmap}
\end{figure*}

 Note that nodes belong to graph $G$ are the actual geographical points, whereas partial nodes belong to $G'$ represent logistics requests need to be served. Graph $G$ is a connected and directed graph which means there exists a path from any point to any other point in the graph $G$. However, the adjacency matrix of graph $G'$ is an upper triangular matrix, because the feasible sequence of serving requests is constrained by the given pickup time $t_i$. \textcolor{black}{ That is, as for an electric vehicle, a request with early pickup time will be served before a later one.} We denote $\mathcal{E}'$ as $\left\{(i, j) | i, j \in \mathcal{V}^{\prime}\right\}$ such that the distance from node $i$ to $j$ can be calculated by the classical shortest path algorithm like Dijkstra's algorithm\cite{wiki_shortest}.  For example we choose node $v_1=16$ and $v_n=68$ as the starting depot and end depot respectively. Node set $v_\mathcal{R}=\{19, 29, 36, 45, 49, 50, 52, 53	\}$ is chosen as transportation request set $\mathcal{R}$. Then a directed acyclic graph, $G'$, is constructed. The relation between graph $G$ and $G^{'}$ is shown clearly in  Fig.\ref{roadmap} \textcolor{black}{ and the travel time (distance) satisfies the property of triangle inequalities for both graph and  reduced graph.}

\par We use $G'(\mathcal{V}',\mathcal{E}')$ instead of $G(\mathcal{V},\mathcal{E})$ in the subsequent sections. This can greatly reduce the number of binary decision variable $x^k_{ij}$ by the reduction of considered vertices and simplify joint routing and charging problem of EVs.

\section{ Mixed Integer Programming Model}
 We establish the joint routing and charging problem of multiple EVs based on a network flow model over graph $G^{'}$ as problem 1\cite{jborlin}. As shown in (\ref{obj}), the objective function minimizes the total charging costs, charging time, travel time and EVs usage fees, and  maximizes the revenue of finishing customer requests, with four constant trade-off parameters $\omega_1, \omega_2, \omega_3$ and $\omega_4$. We can denote the total charging costs and charging time by $\sum_{k\in\mathcal{K}}\sum_{i\in\mathcal{V}'} \sum_{j\in\mathcal{V}'} r^k_{i} g_{i} x^k_{i j}$ and $\sum_{k\in\mathcal{K}}\sum_{i\in\mathcal{V}'} \sum_{j\in\mathcal{V}'} r^k_{i} p_{i}x^k_{ij}$ respectively with charging duration per unit energy $g_i$, charging price $p_i$, binary variable for choosing edge $ij$ for EV $k$, $x^k_{ij}$, and charged amount $r^k_i$ at node $i$ of EV $k$. Besides, the overall travel time of all EVs can be expressed as $ \sum_{k\in\mathcal{K}}\sum_{i\in\mathcal{V}'} \sum_{j\in\mathcal{V}'} T_{i j} x^k_{i j}$ with $T_{ij}$ denoting travel time between node $i$ and $j$.
\textcolor{black}{ With $c_i$ denoting a unified cost vector which represents both the EV usage fee $c_{v}$ and the negative serving cost $-\mathscr{M}_i$:
		$c_i=\left\{
		\begin{array}{ll}
			{-\mathscr{M}_i,} & {\text { if } i\in\mathcal{R}} \\ 
			{c_{v},} & {\text { if } i=v_1}
		\end{array}
		\right.
		$,} the revenue and cost of using EVs can be unified as  $\sum_{k\in\mathcal{K}}\sum_{i\in\mathcal{V}'} \sum_{j\in\mathcal{V}'} c_{i} x^k_{i j}$.  Then the objective function is formulated as formula (\ref{obj}).

 \par \textbf{Problem 1: MIP joint routing and charging problem}
 \begin{equation}\label{obj}
\begin{aligned}
\min\limits_{x_{i j}\in\{0,1\}, r_i^k\in \mathbb{R}  }  \sum_{k\in \mathcal{K}} \Bigg\{ \omega_1 &\sum_{i\in\mathcal{V}'} \sum_{j\in\mathcal{V}'} r^k_{i} g_{i} x^k_{i j}+ \omega_2\sum_{i\in\mathcal{V}'} \sum_{j\in\mathcal{V}'} r^k_{i} p_{i}x^k_{ij}\\
  & +\omega_3 \sum_{i\in\mathcal{V}'} \sum_{j\in\mathcal{V}'} T_{i j} x^k_{i j}\\
  &+\omega_{4}\sum_{i\in\mathcal{V}'} \sum_{j\in\mathcal{V}'} c_{i} x^k_{i j} \Bigg\}  
  \end{aligned}
\end{equation}

  \begin{equation}\label{cons_flow}
\begin{aligned}
\text{s.t.\qquad}  \sum_{j\in\mathcal{V}'} & x^k_{i j}-\sum_{j\in\mathcal{V}'} x^k_{j i}=b_{i}, \quad \forall  i \in \mathcal{V}'; k\in\mathcal{K}\\ &b_{1}=1, b_{n}=-1, b_{i}=0,
\end{aligned}
\end{equation}
\begin{equation}\label{cons_visited}
    \sum_{k\in\mathcal{K}}\sum_{j\in\mathcal{V}'} x^k_{ij} \leq 1,
    \quad\forall i\in\mathcal{R}
\end{equation}
\begin{equation}\label{cons_time}
t_{j} \geq\left(T_{ij}+g_ir^k_{i}+t_{i}\right) x^k_{ij} \quad \forall i\in \mathcal{V}' \setminus v_n,j \in \mathcal{V}' \setminus v_1,k\in\mathcal{K}
\end{equation}
\begin{equation}\label{cons_soc}
E^k_{j}=\sum_{i\in\mathcal{V}'}\left(E^k_{i}+r_i^k-e_{i j}\right) x^k_{i j}, \quad\forall  j\in\mathcal{V}' \setminus v_1,k\in\mathcal{K}
\end{equation}
\begin{equation}\label{cons_socbound}
    0\leq E^k_i\leq E_{k,max}, \quad  i\in \mathcal{V}',k\in\mathcal{K}
\end{equation}
\begin{equation}\label{cons_energybound}
  0\leq  r^k_i\leq r_{max}, \quad  i\in \mathcal{V}', k\in\mathcal{K}
\end{equation}

\par To meet the physical requirements of charging system and  transportation system, network flow conservation constraints, battery energy constraints and pickup time constraints are defined as (\ref{cons_flow})-(\ref{cons_energybound}) in Problem 1.
\begin{itemize}
    \item  \textcolor{black}{Constraints (\ref{cons_flow}) indicates that all EVs are subject to the flow conservation constraint. In other words, EV entering in the request node has to exit out at the same request node, and EV starts at the start depot while returning back to the end depot. In addition, both the single depot vehicle routing problem and the multi-depot vehicle routing problem all can be characterized by (\ref{cons_flow}). The later one can be incorporated into single depot case by adding two virtual depots connecting the start and end depot separately via zero travel distance.}
    \item  Besides, constraint (\ref{cons_visited}) means each customer can be served at most once by the fleet. \textcolor{black}{Note that some transportation
requests may not be served by the fleet if the request associates positive weighted cost and is not beneficial for EV operator.}
    \item Time constraints with prescribed pickup time $t_j$ of customer and charging time $g_i r_i^k$  of EV are specified in (\ref{cons_time}), which states that the arrival time at request nodes should not be later than prescribed pickup time for every transportation request.
    \item Battery energy dynamics of EVs is characterized by equation (\ref{cons_soc}) in which $e_{ij}$ denotes the amount of energy consumption on edges $ij$. $e_{ij}=\epsilon d_{ij}$ where $\epsilon$ is the energy consumption rate and $d_{ij}$ is the travel distance from node $i$ to $j$, $i,j\in \mathcal{V}'$.
    \item Constraints (\ref{cons_socbound}) and (\ref{cons_energybound}) enforce lower and upper bounds on battery energy level of EV $k$, $E_i^k$ and charging amount at node $i$, $r_i^k$. The charging amount at end depot $v_n$, i.e. $r^k_{v_n},\forall k\in\mathcal{K}$ is set as 0. 

\end{itemize}


\par Clearly, Problem 1 is an MIP problem, which can be solved by commercial solvers, such as CPLEX or Gurobi. However, MIP problem is NP-hard and computationally very demanding. In addition, the bilinear term $x_{ij}^k r_i^k$ couples the binary route decision and the continuous charging decision in both objective function (\ref{obj}) and constraints (\ref{cons_time}) and (\ref{cons_soc}), thus making the problem more challenging. To this end, next we propose a two-stage algorithm, which can attain a near-optimal solution quickly, by transforming the original MIP problem with bilinear term to two LP problems. \textcolor{black}{ The main idea behind this transformation lies at a near-optimal solution is attained by optimizing the routing related binary variables and charging related continuous variables separately, which can be further reformulated as two LP problems.}

\section{ Two-stage optimization scheme}
\begin{figure}[ht]
\centering
\includegraphics[width=.95\linewidth]{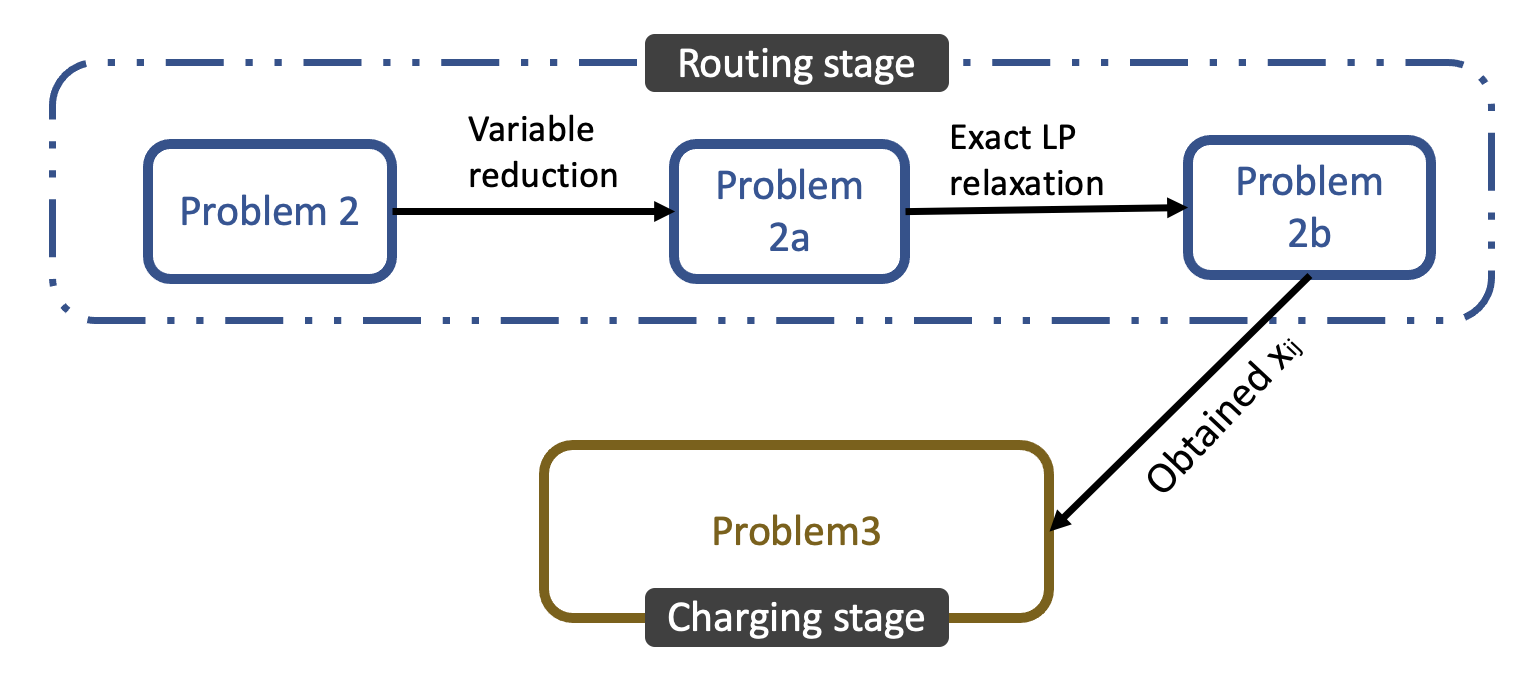}
\caption{\color{black}{Transformation procedure of Problem 2 and 3}}
\label{transchart}
\end{figure}

The proposed two-stage scheme consists of: 1) a routing stage determining the route decision $x_{ij}$, which optimizes travel time, EVs usage fee and service revenue; and 2) a charging stage determining the energy charging amount $r_i$, which optimizes charging time and charging cost. The details are given below. \textcolor{black}{The relationship between the optimization problems of the two-stage optimization scheme is illustrated in Fig.\ref{transchart}.}

\subsection{Routing Stage} 


At this stage, we shall only focus on the routing related objective terms and constraints, while drop charging related terms. That is, in Problem 2 the objective function (\ref{obj_routing}) partially contains the traveling distance and service revenue, and the constraints only include (\ref{cons_flow}-\ref{cons_time}) that are routing related. 
\par \textbf{Problem 2: MIP routing problem}

\begin{equation}
\min\limits_{x^k_{i j}\in\{0,1\}, r_i^k\in \mathbb{R}  } \quad \omega_3 \sum_{k\in\mathcal{K}} \sum_{i\in\mathcal{V}'} \sum_{j\in\mathcal{V}'} T_{i j} x^k_{i j}+\omega_{4}\sum_{k\in\mathcal{K}}\sum_{i\in\mathcal{V}'} \sum_{j\in\mathcal{V}'} c_{i} x^k_{i j}
\end{equation}

$$\text{s.t.\qquad (\ref{cons_flow}-\ref{cons_time})}$$





\par Since all EVs station at the starting depot and return back to end depot after finishing transportation requests, EVs can be treated as a bundled vehicle flow. For the inconvenience of  a three-index formulation in Problem 2, we use two-index vehicle flow model\cite{schneider2014electric} to reformulate Problem 2 to Problem 2a as below, which reduces the number of required variables.

\par\textbf{ Problem 2a: Reformulated MIP routing problem}

\begin{equation}\label{obj_routing}
\min\limits_{x_{i j}\in\{0,1\}, r_i\in \mathbb{R}  } \quad \omega_3 \sum_{i\in\mathcal{V}'} \sum_{j\in\mathcal{V}'} T_{i j} x_{i j}+\omega_{4}\sum_{i\in\mathcal{V}'} \sum_{j\in\mathcal{V}'} c_{i} x_{i j} \end{equation}

 \begin{equation}\label{cons_twoindex_flow}
\begin{aligned}
\text{s.t.\qquad}   \sum_{j\in\mathcal{V}'} & x_{i j}-\sum_{j\in\mathcal{V}'} x_{j i}=b_{i}, \quad \forall  i \in \mathcal{V}'\\ 
&b_{1}=|\mathcal{K}|, b_{n}=-|\mathcal{K}|, b_{i}=0,
\end{aligned}
\end{equation}

\begin{equation}\label{cons_twoindex_visited}
\sum_{j \in \mathcal{V}'} x_{i j}\leq1, \quad \forall i \in \mathcal{R}
\end{equation}

\begin{equation}\label{cons_twoindex_time}
t_{j} \geq\left(T_{ij}+g_ir_{i}+t_{i}\right) x_{ij} \quad \forall i\in \mathcal{V}' \setminus v_n,j \in \mathcal{V}'\setminus v_1
\end{equation}
It can be seen that (\ref{cons_flow}) has been replaced by (\ref{cons_twoindex_flow}) in Problem 2a. Decision variables associated with vertices can be used to keep track of vehicle states and the vehicle index $k$ can be dropped in the Problem 2a\cite{schneider2014electric}.

\par It is clear that Problem 2a is still an MIP problem with a bilinear term $r_i x_{ij}$ in constraint (\ref{cons_twoindex_time}). In order to devise a fast algorithm, we shall first apply exact LP relaxation and then eliminate this bilinear term as detailed below.

\subsubsection{\textit{\textbf{LP relaxation}}}
By relaxing the binary decision $x_{ij} \in \{0,1\}$ to continuous decision $x_{ij} \in [0,1]$, the MIP problem is converted into a continuous optimization problem. Meanwhile, the exactness of LP relaxation is stated below. 

\begin{lemma}
 The LP relaxation of Problem 2a is exact, i.e., the solution of the relaxed problem 2a is equivalent to the solution of the original MIP problem. 
\end{lemma}

\begin{proof}
 As we can see from Problem 2a, formulae (\ref{obj_routing})-(\ref{cons_twoindex_time}) are in the standard form of the minimum cost flow problem\cite{mcfwiki} for which the LP relaxation is exact. As for the addition constraint (\ref{cons_twoindex_time}), we are able to prove that it does not break the LP relaxation exactness. 

\par  Note that $t_i,i\in\mathcal{R} $ is a parameter given by customers rather than the decision variable. Let us reformulate constraint (\ref{cons_twoindex_time}) as follows
 
 $$0\geq \left(T_{ij}+g_ir_{i}-\left(t_{j}-t_{i}\right)\right)x_{ij}$$
 
Simple observations on the sign of $(T_{ij}+g_ir_i-(t_j-t_i))$ can be stated below,
\begin{itemize}
    \item If $\left(T_{ij}+g_ir_{i}-\left(t_{j}-t_{i}\right)\right)> 0$, then $x_{ij}$ must equal to 0 to satisfy constraints (\ref{cons_twoindex_time}) which means that we can remove constraints (\ref{cons_twoindex_time}) when $\left(T_{ij}+g_ir_{i}-\left(t_{j}-t_{i}\right)\right)\geq 0$.
    
    \item If $\left(T_{ij}+g_ir_{i}-\left(t_{j}-t_{i}\right)\right)\leq 0$, constraints (\ref{cons_twoindex_time}) will be satisfied no matter what value $x_{ij}$ is, thus the constraint (\ref{cons_twoindex_time}) is a redundant constraint.
\end{itemize}
\textcolor{black}{In a word, no matter what value $(T_{ij}+g_ir_i-(t_j-t_i))$ is,  the value of $x$ is restricted  as a binary value.} The exactness of LP relaxation of problem 2a is maintained and the decision variables $x_{ij}$ can be relaxed as the continuous variables. This completes the proof.
\end{proof}

\subsubsection{\textbf{Bilinear term}} After LP relaxation, there is still a bilinear term $r_ix_{ij}$ in constraint (\ref{cons_twoindex_time}). In order to ease the computation burden, $r_i$ is replaced by $e_{ij}$ which eliminates the coupled bilinear. $e_{ij}$ is the energy needed from node $i$ to node $j$. Thus, it implies that we always charge enough energy in order to arrive at the next request node. Then constraint (\ref{cons_twoindex_time}) is reformulated as formula (\ref{cons_linear_time}).
\begin{equation}\label{cons_linear_time}
t_{j} \geq\left(T_{ij}+g_ie_{ij}+t_{i}\right) x_{ij} \quad \forall i\in \mathcal{V}' \setminus v_n,j \in \mathcal{V}' \setminus v_1
\end{equation}

Finally, we obtain the following LP problem without bilinear term for the routing stage.
\par \textbf{Problem 2b: LP routing problem}
\begin{equation*}
\begin{array}{cc}
      \min\limits_{x_{i j}\in[0,1]  }  \qquad & \omega_3 \sum_{i\in\mathcal{V}'} \sum_{j\in\mathcal{V}'} T_{i j} x_{i j}+\omega_{4}\sum_{i\in\mathcal{V}'} \sum_{j\in\mathcal{V}'} c_{i} x_{i j} \\
     \text{s.t.} &(10)-(11) \quad\text{and}\quad (13)
\end{array}
\end{equation*}

  \subsection{Charging Stage}  
  
  After solving Problem 2b by any efficient algorithm, e.g., simplex method, denote the optimal routes by $P=\{P_1,\ldots,P_K\}$, where $P_k = \{v_1,R_k,v_n\}$ is the specific route plan for EV $k,k\in\mathcal{K}$, and $R_k \subseteq R$ is the requests assigned to the $k^{th}$ EV. Conditioned on the optimal route plan, the charging stage optimizes the amount of charged energy, which minimizes the charging time and cost, as formulated in Problem 3.

\textbf{Problem 3: LP charging problem}
 \begin{equation}
   \min\limits_{r_{i}\in\mathbb{R},  i \in P_k}   \omega_1 \sum_{i\in P_k} r_{i} g_{i} + \omega_2 \sum_{i\in P_k} r_{i} p_{i}  
 \end{equation}

\begin{equation}\label{cons_charging_soc}
\text{s.t.}\quad E_{j}=E_{i}+r_{i}-e_{i j}, \forall  j\in P_k \setminus v_1
\end{equation}

\begin{equation}\label{cons_charging_socbound}
{0 \leq E_{i} \leq E_{k,max},}  {\quad  \forall i \in P_k}    
\end{equation}

\begin{equation}\label{cons_charging_energybound}
{0\leq r_{i} \leq e_{ij} \quad  \forall i \in P_k \setminus v_n}     
\end{equation}

\noindent Battery and charging constraints (\ref{cons_charging_soc}-\ref{cons_charging_energybound}) are analogous to (\ref{cons_soc}-\ref{cons_energybound}) with a specified route $P_k$ for each EV. It is clear that Problem 3 is again an LP problem and can be quickly solved.



To summarize, the two-stage algorithm for EVs routing and charging is described below \textcolor{black}{by Algorithm 1, and a flowchart is shown in Fig. 2 in the next section.} With the proposed two-stage optimization scheme, the solution of EVs routing and charging problem can be obtained  within polynomial time. It is notable that the solution is near-optimal because we, 1) regard the original joint optimization problem as two sequential subproblems, and 2) replace $r_j$ by $e_{ij}$ for every request to eliminate the bilinear term. However, in the next section we demonstrate that not every $r_j$ needs to be replaced by $e_{ij}$, and an iterative algorithm is proposed to improve the solution optimality of two-stage method.

\begin{algorithm}[ht]
 \caption{Two-stage optimization algorithm}
 \begin{algorithmic}[1]
 \renewcommand{\algorithmicrequire}{\textbf{Input:}}
 \renewcommand{\algorithmicensure}{\textbf{Output:}}
  \STATE Solve Problem 2b  
  \FOR{$k\in\mathcal{K}$}
  \STATE Solve Problem 3 for each $k
 \in\mathcal{K}$
  \ENDFOR
 \end{algorithmic} 
 \end{algorithm}

\section{ Iterative optimization scheme }

\par Replacing $r_i$ by $e_{ij}$ in (\ref{cons_linear_time}) implies we always charge enough energy at node $i$ before heading to the next node $j$. However, as in practice an EV has stored sufficient energy for most journeys, the replacement 
for all requests is not necessary but too conservative, consequently lowers the quality of solution. To tackle this issue, an iterative variant algorithm of the two-stage one is proposed, which only replaces $r_i$ by $e_{ij}$ for a portion of nodes at each iteration, and yields a much better solution.
\par We firstly define the variable $\Phi^k_j$ indicating whether the actual pickup time of node $j$ satisfies predefined pickup time as below. 

\begin{equation}\label{Set_violated_time}
  \Phi^k_j=\left(T_{ij}+g_ir_{i}+t_{i}-t_{j}\right) x_{ij} \quad \forall i\in \mathcal{V}' \setminus v_1 ,j \in \mathcal{R}  
\end{equation}

It is clear that $\Phi_j^k >0$ implies $t_j <T_{ij}+g_ir_i+t_i$, i.e., the pickup time of node $j$ is too early and cannot be met by EV $k$ after charging and leaving from node $i$. We call such request $j$ as \textit{violated request},

\begin{equation}\label{Set_violated_node}
\mathcal{R}^q_v=\{j|\Phi^k_j > 0,k\in\mathcal{K},j\in\mathcal{R} \}\end{equation}

\noindent where $q$ is the index of iteration. On the other hand, the request $j$ is called as \textit{Non-violated request} if $\Phi_j^k \leq0$.

Noticing that the optimal $r_i$ always satisfies $0\leq r_i\leq e_{ij}$, replacing $r_i$ by $e_{ij}$ in (\ref{cons_linear_time}) makes the constraint conservative and produces a less optimal solution; on the other hand, replacing $r_i$ by $0$ may yield a better solution but could violate the pickup time constraint sometime. Therefore, the basic idea of the iterative variant method comes as the following.
\begin{itemize}
   \item At the routing stage, we can initially use $0$ to replace $r_i$ for all requests, and find the routes.
   \item Next, we can determine $r_i$ at the charging stage, and then use (\ref{Set_violated_time}) to check if there is pickup time violation for a request. 
   \item If no violation exists, the algorithm has already found a solution.
   \item Otherwise, we shall replace $r_i$ by $e_{ij}$ for those violated requests to eliminate the violation and re-solve the routing problem again.
\end{itemize}
\par
\begin{figure}[ht]
\centering
\includegraphics[width=1\linewidth]{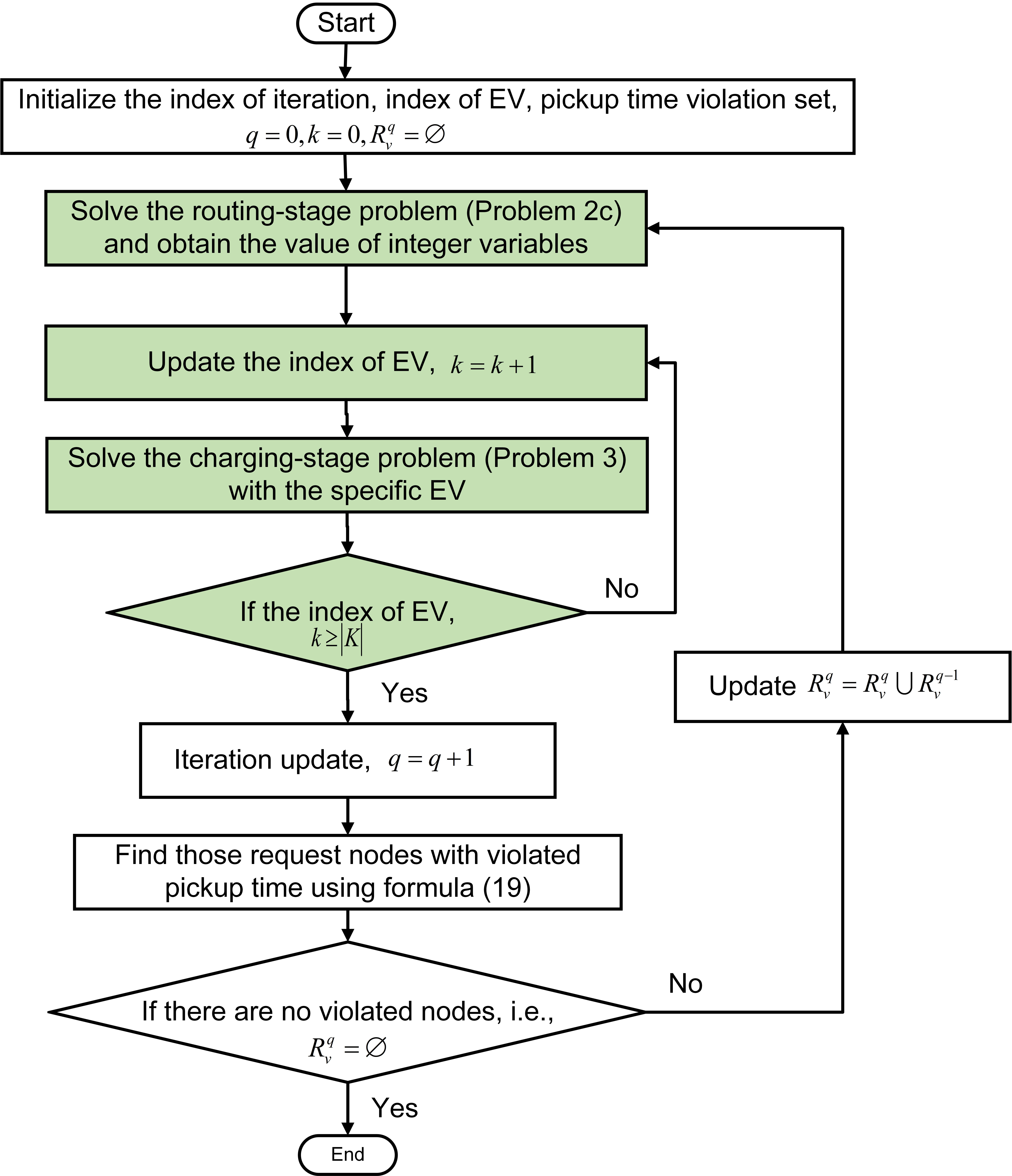}
\caption{\color{black}{Flowchart of algorithms 1 and 2}}
\label{flowchart}
\end{figure}

 In other words, the routing problem in iterative method can be written as below.\\
\textbf{Problem 2c: Iterative LP routing problem}
\begin{equation*}
\begin{array}{cc}
      \min\limits _{x_{i j}\in[0,1]  }   \qquad &\omega_3 \sum_{i\in\mathcal{V}'} \sum_{j\in\mathcal{V}'} T_{i j} x_{i j}+\omega_{4}\sum_{i\in\mathcal{V}'} \sum_{j\in\mathcal{V}'} c_{i} x_{i j}\\
     \text{s.t.} &(10)-(11) \quad\text{and}
\end{array}
\end{equation*}

     \begin{equation}\tag{20-a}
t_{j} \geq\left(T_{ij}+t_{i}\right) x_{ij} \quad \forall i\in\mathcal{R}, j \in \mathcal{R} \setminus \mathcal{R}^q_v
 \end{equation}
 \begin{equation}\tag{20-b}
t_{j} \geq\left(T_{ij}+g_ie_{ij}+t_{i}\right) x_{ij} \quad \forall i \in \mathcal{R},j\in\mathcal{R}^q_v
 \end{equation}
Above iterative process will stop when no pickup time violation exists, i.e., $R_v^q = \emptyset$ for the $q^{th}$ iteration. Detailed procedures of the iterative optimization scheme is
shown in Algorithm 2 and the flowchart of Algorithm 2 is provided in Fig.\ref{flowchart}. \textcolor{black}{The green part can describe Algorithm 1 except that Problem 2c is replaced with Problem 2b in Algorithm 1.}



\begin{algorithm}[ht]
 \caption{Iterative optimization algorithm}
 \begin{algorithmic}[1]
 \renewcommand{\algorithmicrequire}{\textbf{Input:}}
 \renewcommand{\algorithmicensure}{\textbf{Output:}}
 \STATE  
 \textbf{\textit{Initialization} }:$q=1;\mathcal{R}^q_v=\emptyset;$ 
  \STATE Solve Problem 2c  
  \FOR{$k\in\mathcal{K}$}
  \STATE Solve Problem 3 for each $k
 \in\mathcal{K}$
  \ENDFOR
   \STATE \color{black}{$q = q+1$} \color{black}
  \STATE \textit{Violation check}: calculate  $\mathcal{R}^q_v=\{j|\Phi^k_j > 0,k\in\mathcal{K},j\in\mathcal{R} \}$
  
  \IF{ $R_v^q = \emptyset$}
    \RETURN optimization result
    \ELSE
\STATE  Update $R_v^q = R_v^q \cup R_v^{q-1}$

\STATE Go to Step 2
\ENDIF
  
    


 \end{algorithmic} 
 \end{algorithm}

    




\begin{lemma}
Algorithm 2 converges in finite iterations.
\end{lemma}

\begin{proof}
 The idea behind iterative optimization scheme lies in adding time constraints (20-b) on request nodes with violated pickup time iteratively. Then the iterative scheme will certainly terminate within finite steps because of finite number of transportation requests $\mathcal{R}$. 
\end{proof}

\par   Since Problem 2b, Problem 2c and Problem 3 are all LP problems and can be solved in polynomial time, the computation speed of two-stage optimization scheme for EVs routing and charging problem is much faster than solving MIP problem. The computation time of iterative optimization scheme is slightly longer than the two-stage one with the improvement in solution quality as illustrated later.

\begin{remark}\label{remark_2}
It is guaranteed that the iterative LP method (Problem 2c) can obtain a solution no worse than the LP method (Problem 2b) at the routing stage. It is seen that Problem 2c is similar to Problem 2b, except for replacing constraint (\ref{cons_linear_time}) by (20). Since (20) is less tight than (\ref{cons_linear_time}), Problem 2c has a larger feasible region. Therefore, the solution of Problem 2c  is certainly no worse than that of Problem 2b. Although we expect that a better solution at routing stage could lead to a better solution for the entire two-stage sequential decision problem, it is not guaranteed and a few exceptions are observed in the simulation. In fact, it is possible that a worse charging solution is obtained at stage 2, even under a better routing solution of stage 1.
\end{remark}



\section{Numerical Simulation}
 In this section, six simulation scenarios, varying in terms of  the size of transportation network, costs of EV usage, are employed to carry out EVs routing and charging simulation to validate the effectiveness of the proposed algorithms. There are three transportation networks considered as connected graphs for simplicity. All of them have 100 road intersections. For ease of notation, \textit{ILP, TLP and MIP} represent solving joint routing and charging problem of EVs with iterative LP, two-stage LP and MIP method, respectively.
 
 \subsection{Parameter Settings}
 The common parameter settings are introduced firstly. There are two kinds of charging nodes, one with fast charging duration per unit energy $g_2$ and another one with slow charging duration per unit energy $g_1$ where $g_1=\frac{36.67}{60} \text{min/kWh} $ and $ g_2=\frac{1}{12} \text{min/kWh}$. The charging price for all charging nodes is set as $\$\frac{1}{12} $ \cite{Pourazarm}.  All four weight factors $\omega_1,\omega_2,\omega_3,\omega_4$ are fixed as 1. Although the results can be easily extended to heterogeneous scenario, we assume all EVs are homogeneous and the energy capacity of EV is $E_{k,max}=100\text{kWh},\forall k\in \mathcal{K}$. Initial energy of all EVs at the depot is the same, i.e. $E_k^0=0.7 E_{k,max},\forall k\in\mathcal{K}$. The costs of EV usage have two different values $\$200 $ and $\$500 $. Revenues of serving transportation requests are assumed as the same, which is  $c_i=\$10 ,i\in\mathcal{R}$. Pickup time of requests are given by customers.

\par Three transportation networks are modified from Desaulniers' EVRP benchmark\cite{desaulniers2016exact}. Instances c106\_21\_50, c101\_21\_50 and rc101\_21\_50 \cite{Errico} are retrieved from webpage \url{https://w1.cirrelt.ca/~errico/#Instances} and modified as three  different transportation networks. We extract  Column \textit{x,y and ReadyTime} with \textit{Row C1-C50} from retrieved instances to form new transportation networks with the value of Column \textit{ReadyTime} assigned as pickup time of transportation requests, $t_i,i\in\mathcal{R}$. The detailed parameter settings of  six scenarios are summarized in Table \ref{dataset}. 


\begin{table*}[ht]
    \centering
    \begin{tabular}{c|c|c|c|c|c|c}
    \hline
    \hline
        Scenario & \#1 &\#2 &\#3 &\#4 &\#5 &\#6  \\
        \hline
        \hline
        Map & c106 & c106 & c101 & c101 & rc101 & rc101  \\  
       EV usage fee (\$)  & 200 & 500& 200 & 500& 200 & 500 \\  
        Revenue of serving requests (\$) &10&10&10&10&10&10  \\  
       Initial SOC of EV &0.7 &0.7&0.7&0.7&0.7&0.7 \\  
       Energy consumption rate (kWh/km)  &0.23 &0.23&0.23&0.23&0.23&0.23\\  
        Speed (km/h)   &18 &18 &18 &18 &18 &18 \\  
        Number of EV &2 &2 &2 &2 &2 &2\\
         \hline
         \hline
         \end{tabular}
    \caption{Simulation data}
    \label{dataset}
\end{table*}



    

\par Several experiments are carried out with three different solution methods, i.e. MIP, ILP and TLP. The performances of proposed methods in terms of solution quality and computational speed are compared with the MIP method, \textcolor{black}{which can obtain a global optimal solution for joint  routing and charging problem of EVs. Commercial solver, CPLEX, is employed to solve this MIP problem to optimality\cite{gay2011ibm}.}   Size of $\mathcal{R}$ is chosen from \{17, 19, 21, 23, 25, 27, 29, 31, 33, 35\} sequentially to investigate the performance under different number of requests,  and results are shown in Fig.\ref{ObjValSce1}-\ref{SerRequSce1}. All optimization methods are implemented with MATLAB R2018b on a PC with Intel Core i5-7500 3.4 GHz and 8 GB of RAM.  

\subsection{Simulation Results}

\par  We use three indices: computation time, solution quality and served requests to evaluate performance of TLP and ILP as shown in Fig.\ref{ObjValSce1}-\ref{SerRequSce1} under parameter settings of scenario 1. Simulation results of all 6 scenarios are summarized in Fig.\ref{Sol6Sce}-\ref{Ser6Sce}.

\subsubsection{Solution quality}
\begin{figure}[ht]
\centering
\includegraphics[width=.65\linewidth]{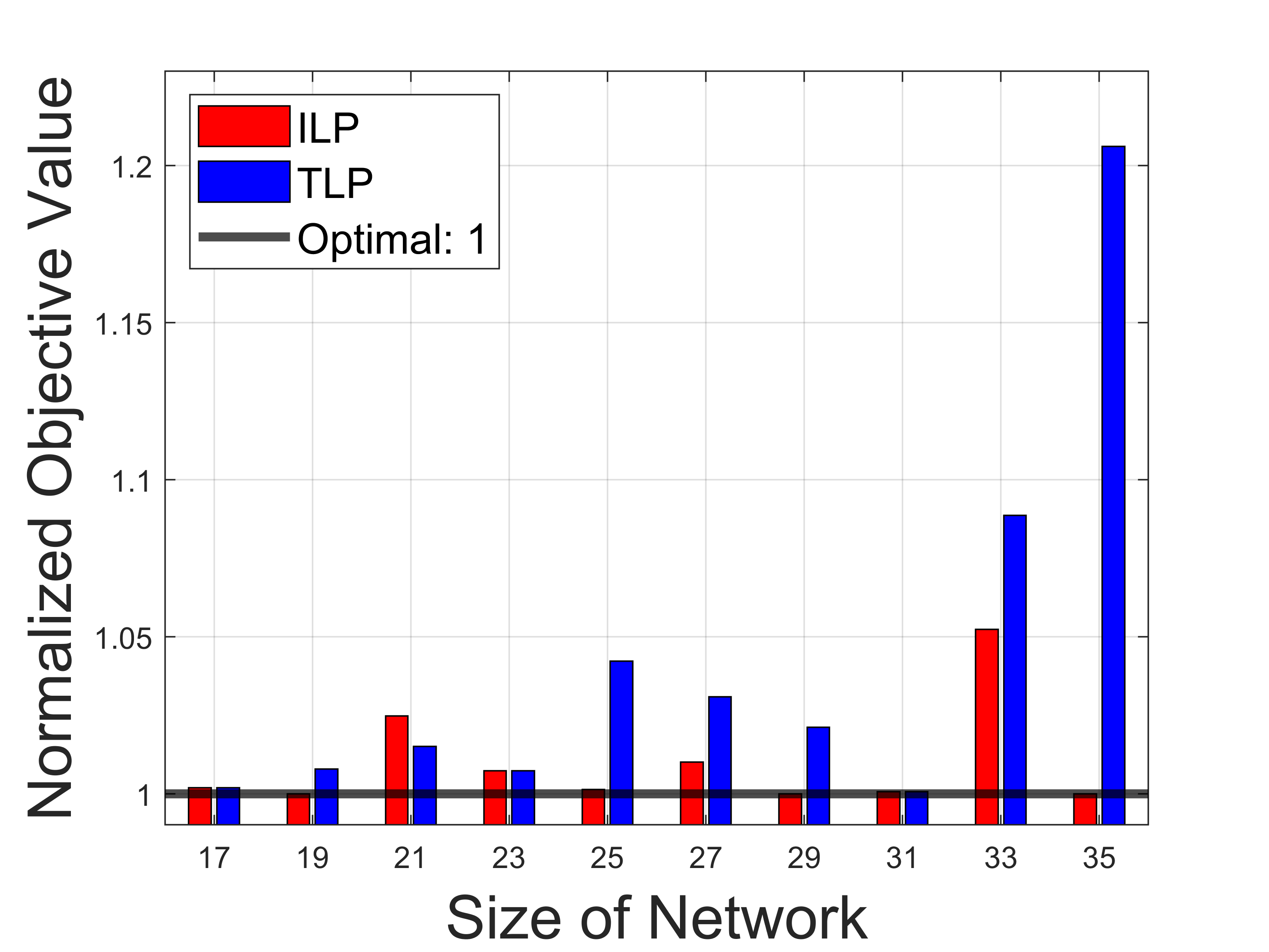}
\caption{Objective value for Scenario 1}
\label{ObjValSce1}
\end{figure}

Solution quality of ILP, TLP and MIP is shown in Fig.\ref{ObjValSce1} Objective values of ILP and TLP are normalized by the optimal solution obtained from MIP. In most cases, solution quality of ILP is better than TLP, especially in large-scale problem. However, exception has been observed when the number of requests is 21, and the reason is explained by Remark \ref{remark_2}. As illustrated in Fig.\ref{ObjValSce1}, the maximum deviation of ILP value from the MIP optimal value is less than 5 percent while the maximum deviation of TLP value can reach 15 percent, which validates that the ILP outperforms TLP in terms of solution quality. 

\subsubsection{Computation speed}
\begin{figure}[ht]
\centering
\includegraphics[width=.65\linewidth]{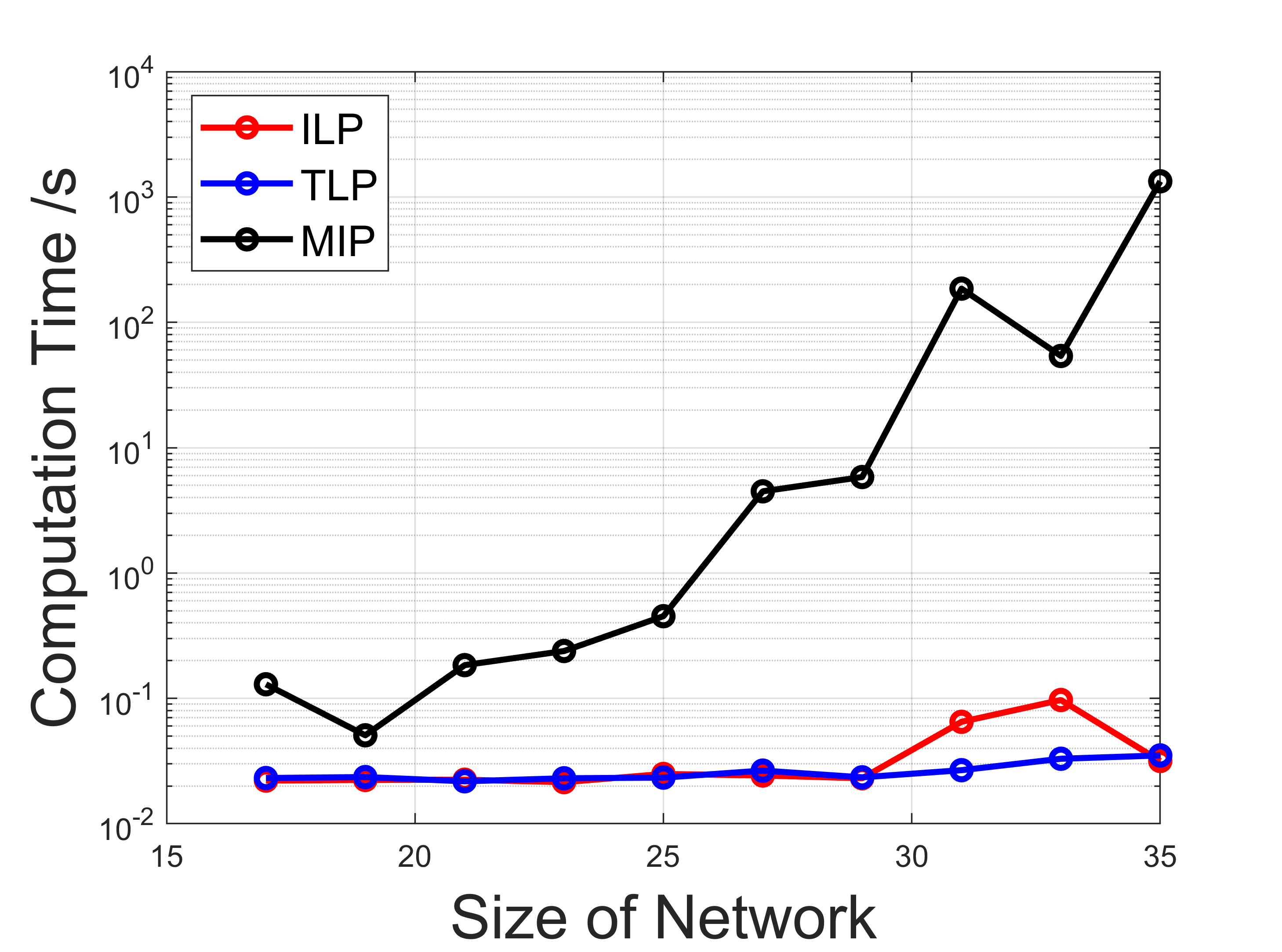}
\caption{Computation time for Scenario 1}

\label{ComTimeSce1}
\end{figure}

Computation time of ILP, TLP and MIP are shown in Fig.\ref{ComTimeSce1}.  The computation time of MIP grows exponentially with the size of requests set $\mathcal{R}$ which makes it difficult to solve the large scale joint routing and charging problem. On the other hand, both TLP and ILP achieve much faster computational speed.


\subsubsection{Served requests}
\begin{figure}[ht]
\centering
\includegraphics[width=.65\linewidth]{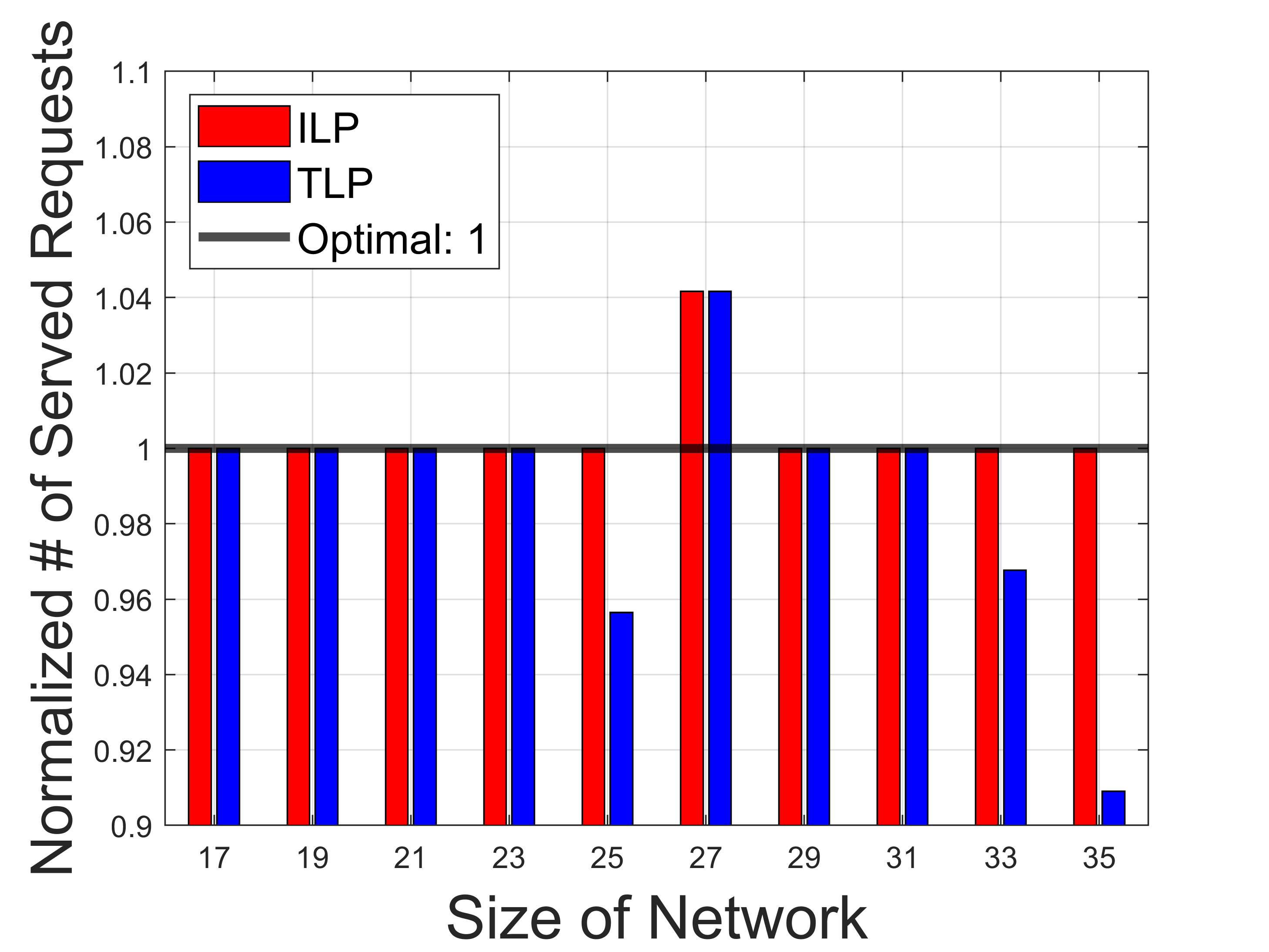}
\caption{Served requests of Scenario 1}
\label{SerRequSce1}
\end{figure}

Served requests of ILP, TLP and MIP are shown in Fig.\ref{SerRequSce1}. The number of served requests for ILP and MIP nearly has the same value for most cases except that sometimes the number of served requests of ILP is slightly larger than MIP as shown in Fig.\ref{SerRequSce1}. This situation can be explained as some EVs drive longer distance in order to serve more transportation requests. However, revenues from serving those requests are less than the cost for extra distance resulting in serving more requests with larger objective value.

\subsubsection{All 6 scenarios} 
\begin{figure}
\centering
\includegraphics[width=.65\linewidth]{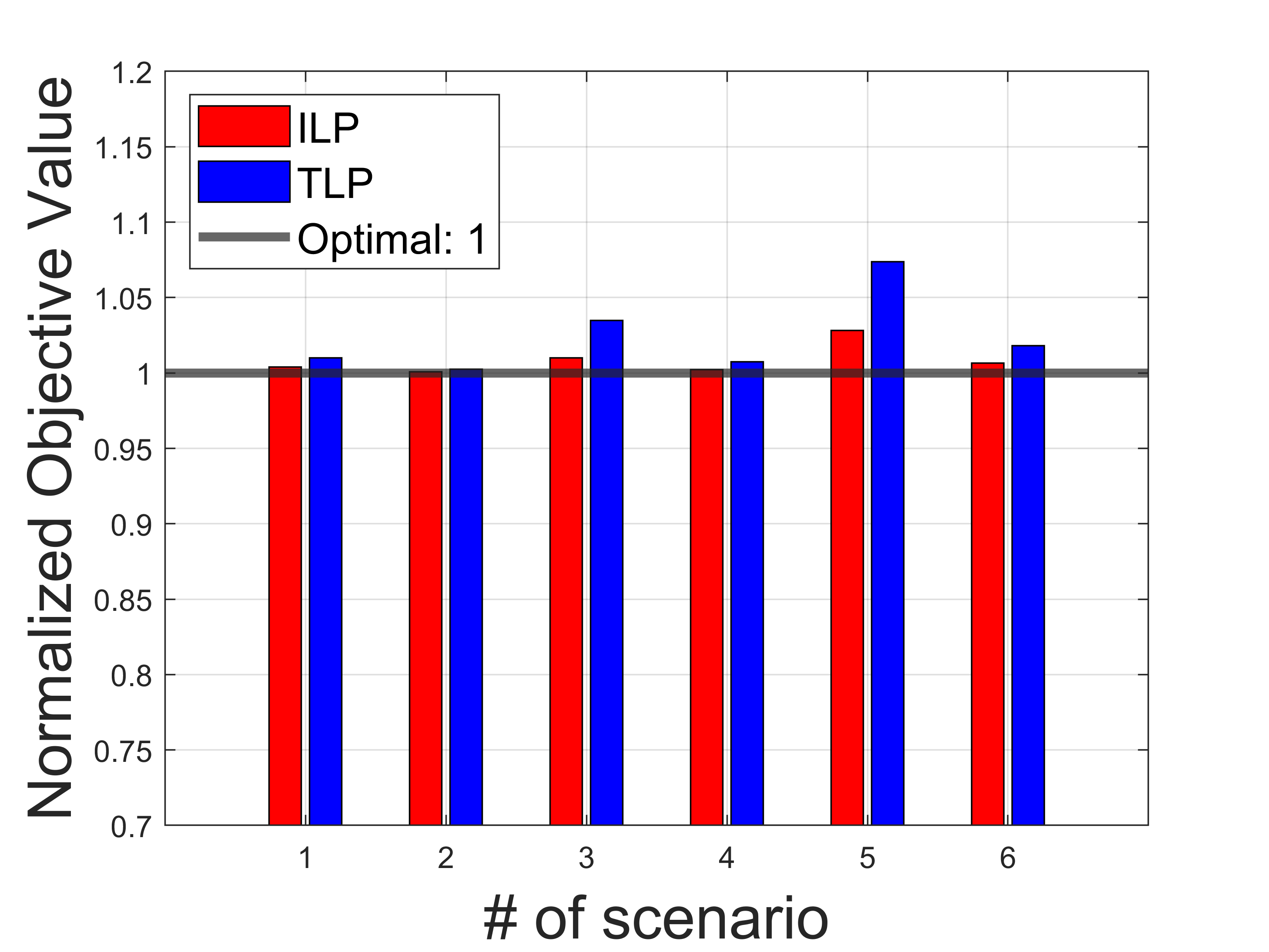}
\caption{Solution quality of 6 scenarios }
\label{Sol6Sce}
\end{figure}

\begin{figure}
\centering
\includegraphics[width=.65\linewidth]{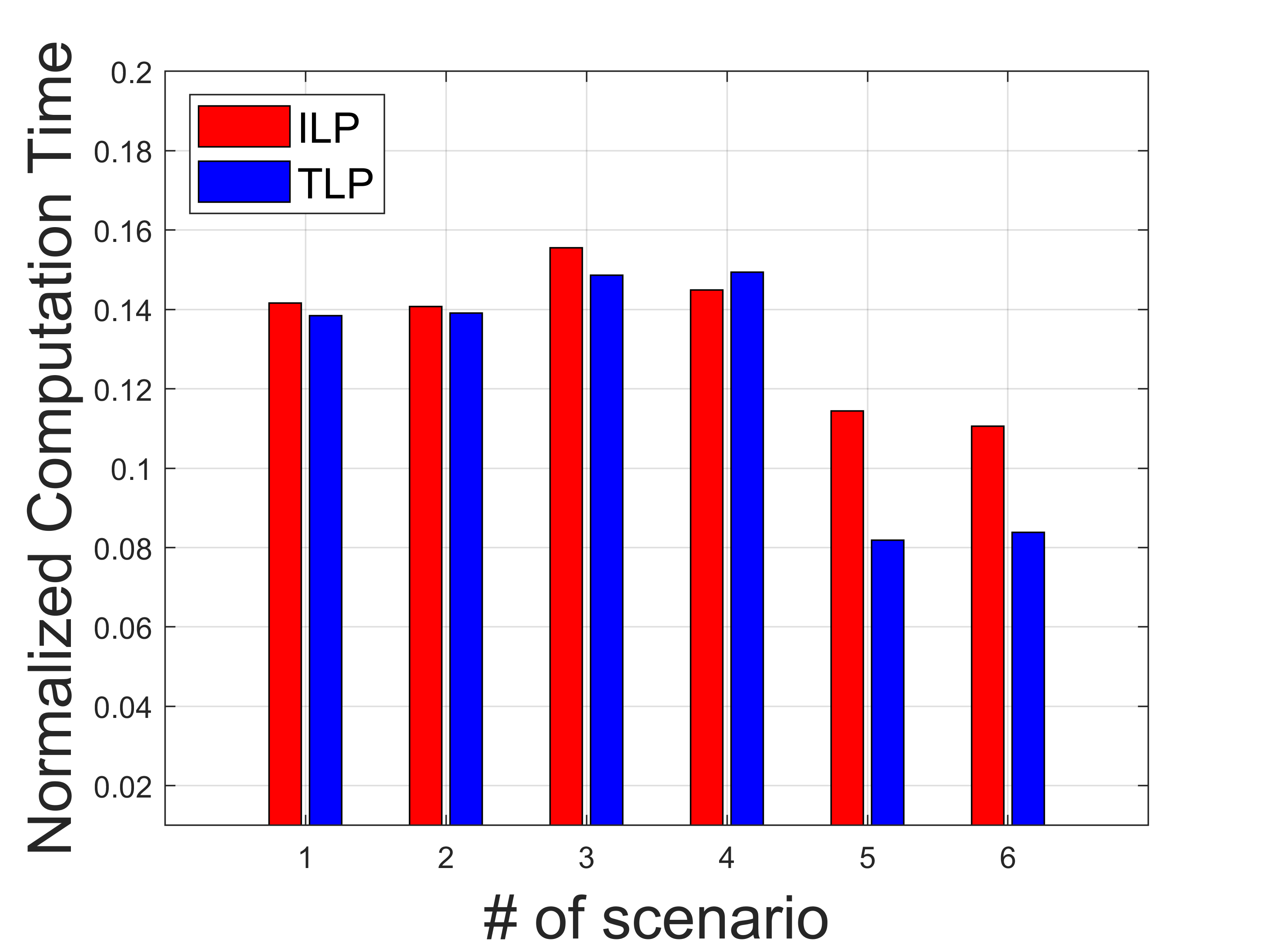}
\caption{Computation speed of 6 scenarios}
\label{Speed6Sce}
\end{figure}

\begin{figure}
\centering
\includegraphics[width=.65\linewidth]{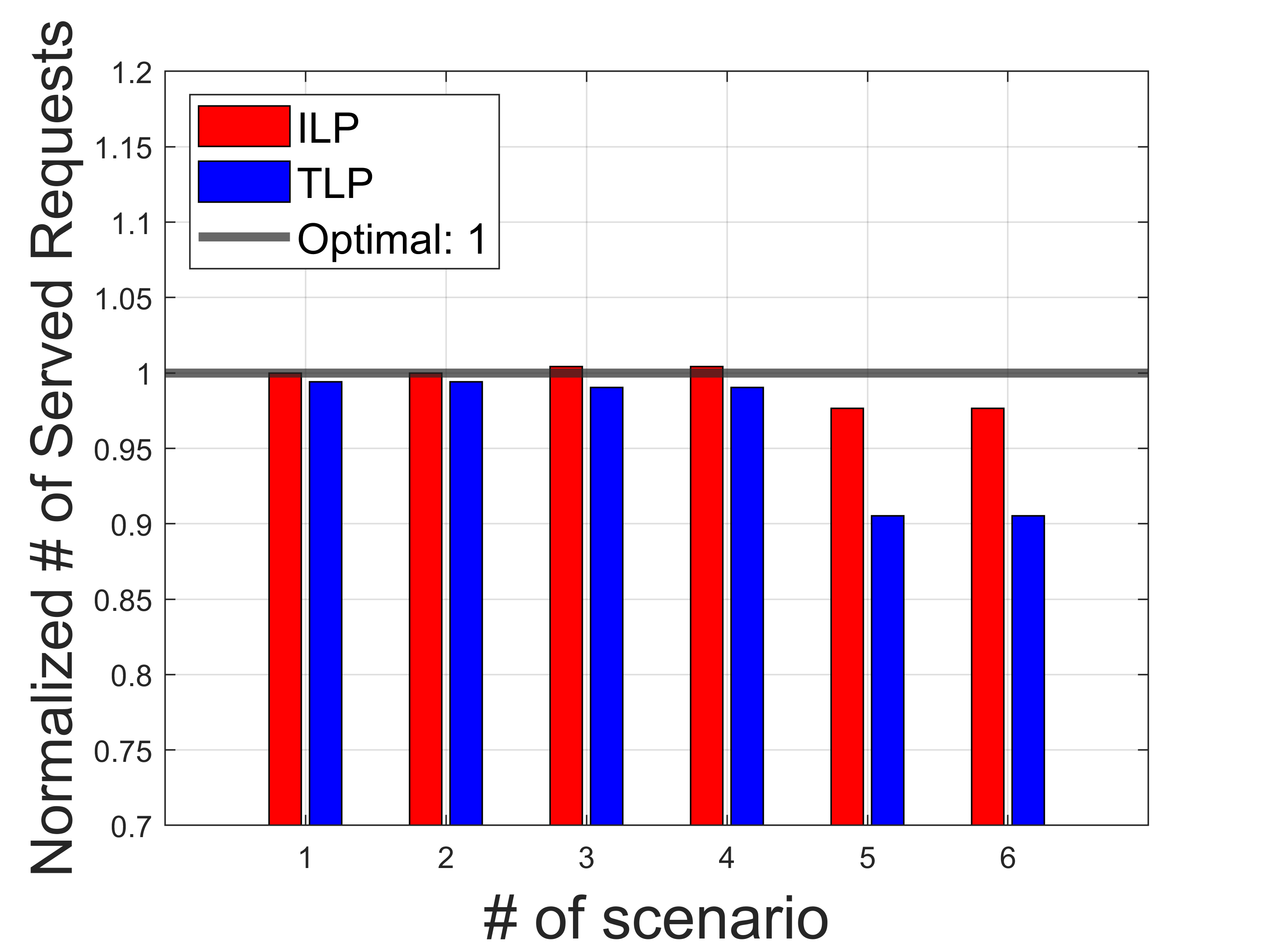}
\caption{Served requests of 6 scenarios}
\label{Ser6Sce}
\end{figure}

Statistical results of 6 scenarios are shown in Fig.\ref{Sol6Sce}-\ref{Ser6Sce}. In these figures, objective values, computation time and the number of served requests of ILP and TLP are normalized by the optimal solution obtained from MIP and then averaged for different size of requests.  As shown in Fig.\ref{Sol6Sce} and Fig.\ref{Ser6Sce},  better performance (solution quality and number of served requests) is achieved by ILP instead of TLP under all 6 scenarios, which validates the effectiveness of ILP. Besides, there is no significant change in computation time for all 6 scenarios with either ILP or TLP method, shown in Fig.\ref{Speed6Sce}.


\subsection{Parameter Sensitivity Test}
In this part, the robustness of proposed algorithms is evaluated under the change of energy consumption rate and traveling speed. The scalability of proposed methods is also discussed. Other parameters are set as Scenario 1 and the number of requests is set as 27.


\par
\begin{figure}
\centering
\includegraphics[width=.65\linewidth]{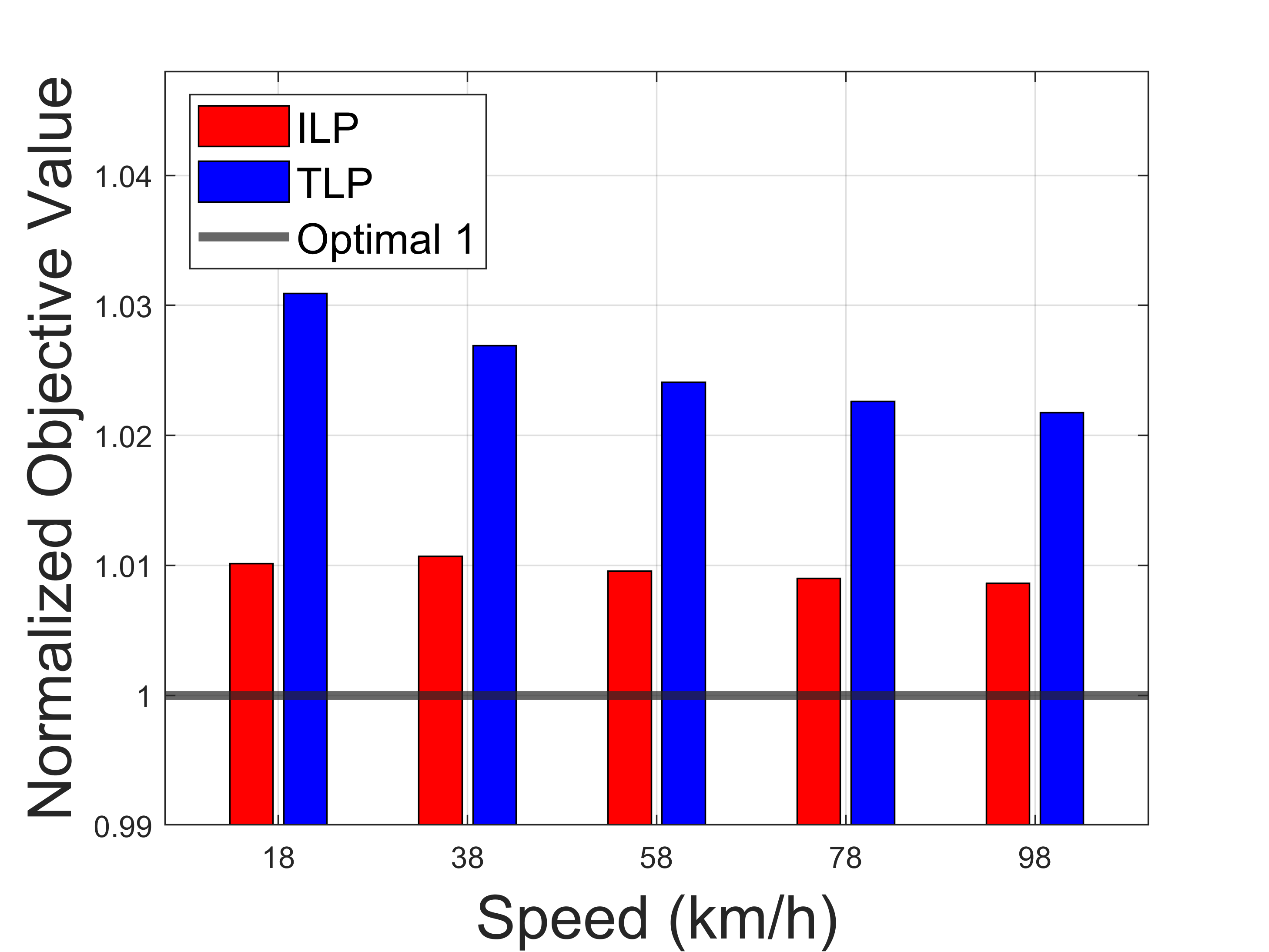}
\caption{Solution quality of ILP and TLP against MIP with Different Speed of EV }
\label{DifferentSpeed}
\end{figure}

 \textbf{ Various Speed:} We implement ILP, TLP and MIP methods with EVs speed changed from 18km/h to 98km/h, representing a reasonable speed range of a logistic EV. As illustrated in Fig. \ref{DifferentSpeed}, ILP and TLP perform very well under all circumstances. The maximum deviations of ILP and TLP from the optimal solution remain within $1.5$ and $3.5$ percent, respectively. 
    
\begin{figure}
\centering
\includegraphics[width=.65\linewidth]{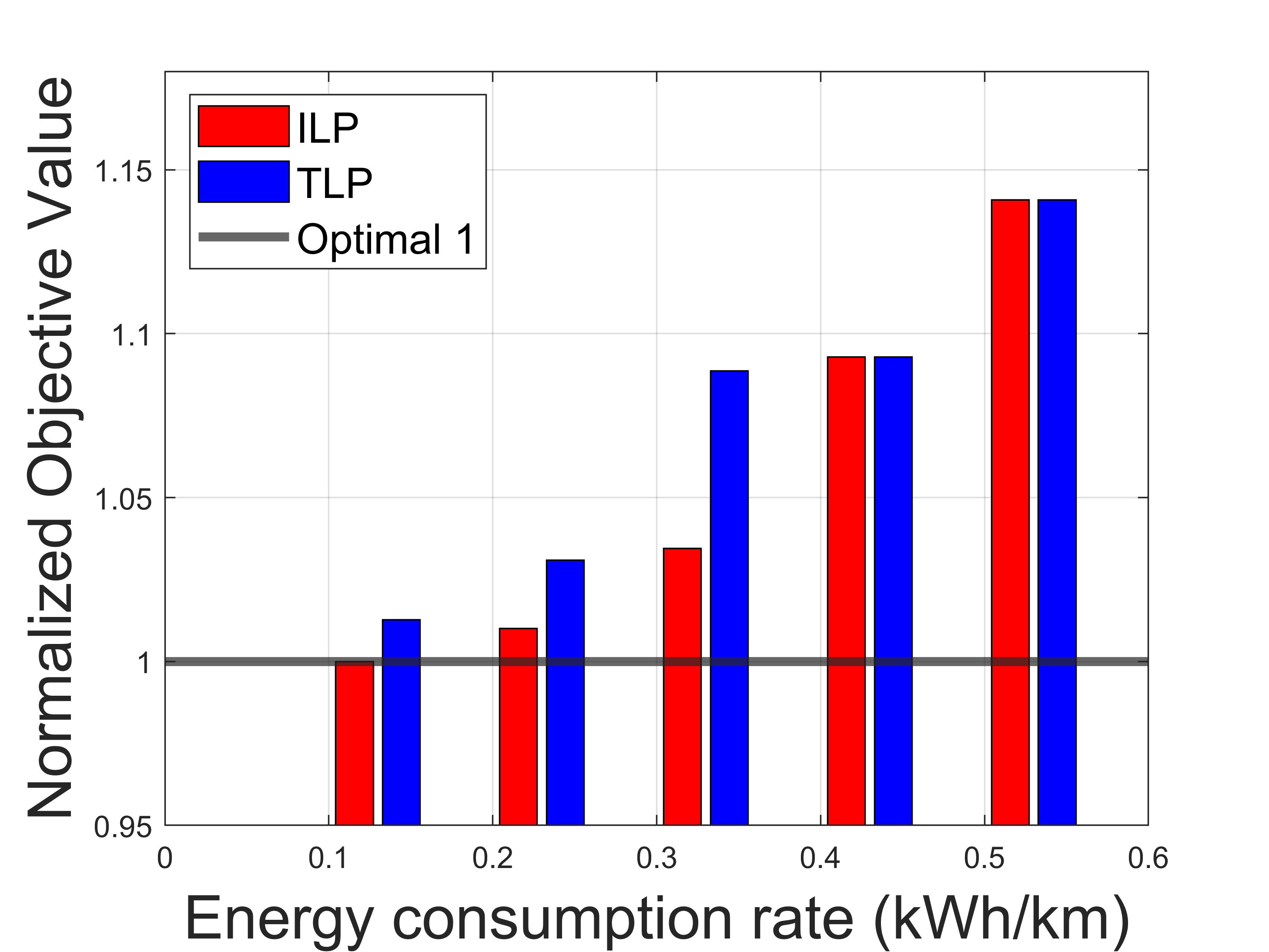}
\caption{Solution quality of ILP and TLP against MIP with Different Energy Consumption Rate of EV}
\label{DifferentEnergy}
\end{figure}

 \textbf{ Various Energy Consumption Rate:} To test the performance of ILP and TLP in terms of energy consumption rate of EV, simulations are carried out with energy consumption rate changing from  0.13 kWh/km  to 0.53 kWh/km, which covers the practical energy consumption range of personal EV, van and light truck\cite{ev_database}.  As illustrated in Fig.\ref{DifferentEnergy}, the performance of ILP and TLP deteriorate with the increase of energy consumption rate. Because large energy consumption rate requires more frequent charging, which renders the proposed methods more conservative.
 
\begin{figure}
\centering
\includegraphics[width=.65\linewidth]{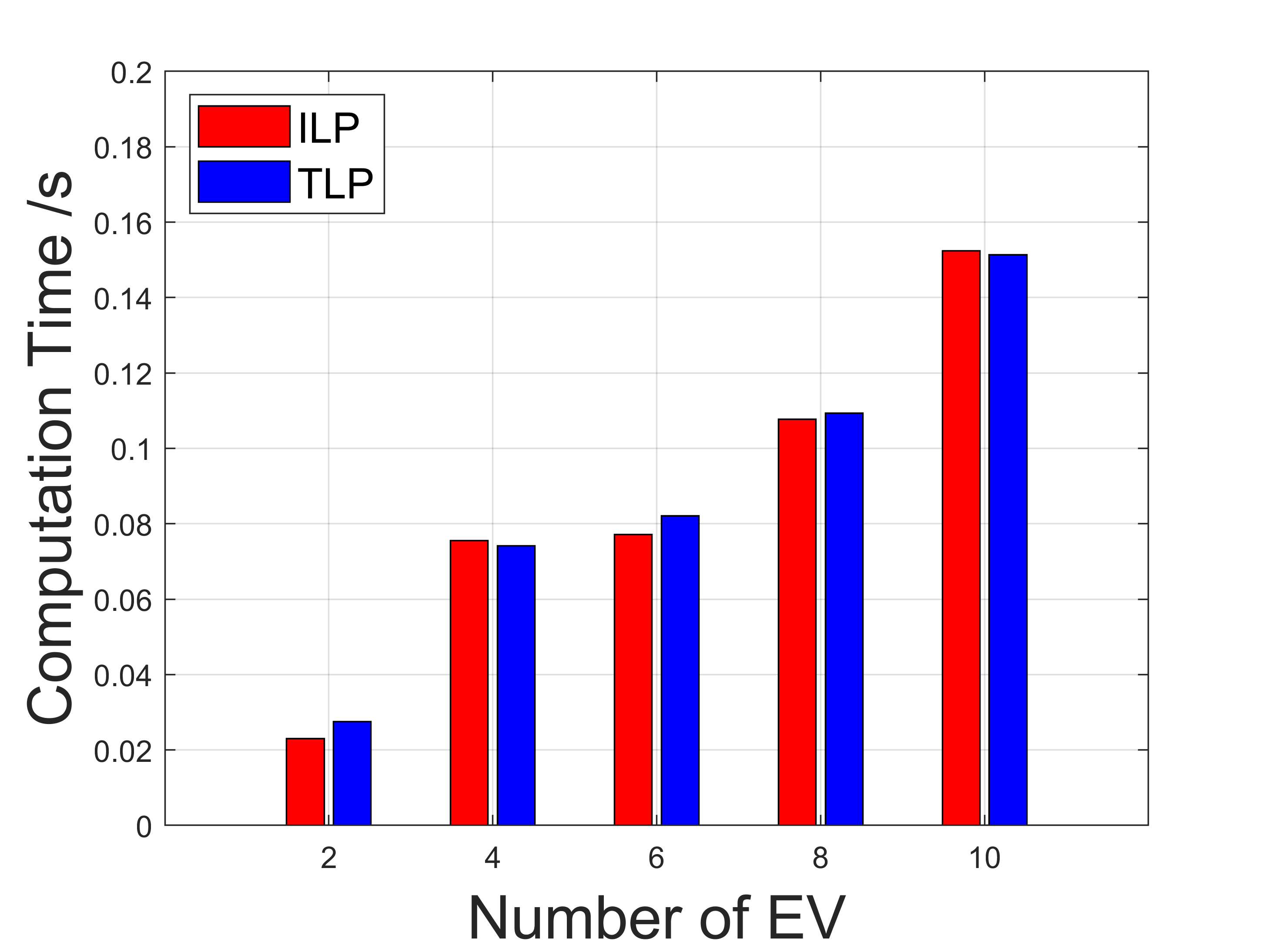}
\caption{Computation speed of ILP and TLP with Different Number of EV}
\label{DifferentNum}
\end{figure}

\textbf{Various Number of EV:} Fig.\ref{DifferentNum} characterizes the scalability of proposed methods in the number of EVs and requests. Computation time of ILP and TLP varies insignificantly with the change of the number of EV from 2 to 10 and the corresponding number of requests increases from 19 to 99, which verify the wide scalability of proposed methods.

\textcolor{black}{\subsection{Case Study on the Road Map of Belgium}
We first adopt the map data of Belgium \cite{BelgiumMap}, which records the geographic coordinates of 50 nodes, and construct its transportation network.	
The battery capacity and charging duration per unit energy are selected as 90 $kWh$, $\frac{1}{22} \text{hr/kWh}$ / $\frac{1}{3} \text{hr/kWh}$ \footnote{Fast and slow charging rates are 22 $kW$ and 3 $kW$ respectively. } (Tesla Model X \cite{ModelX}). 
In addition, the energy consumption rate and average speed of EVs are set as 0.24 $kWh/km$ and 90 $km/h$ \cite{ModelX} respectively. According to the electricity price of Belgium\cite{ElectricityPrice}
,  there are two kinds of charging price $\$$0.299$/kWh$ and $\$$0.129$/kWh$, respectively. The usage cost of EV and revenue from serving a request are set as $\$1299$ and  $\$9.05$ respectively\cite{RequestServe,EVborrow}. Due to the memory limitation to implement CPLEX, the number of considered nodes is restricted under $27$.}
\par
\begin{figure}
\centering
\includegraphics[width=.65\linewidth]{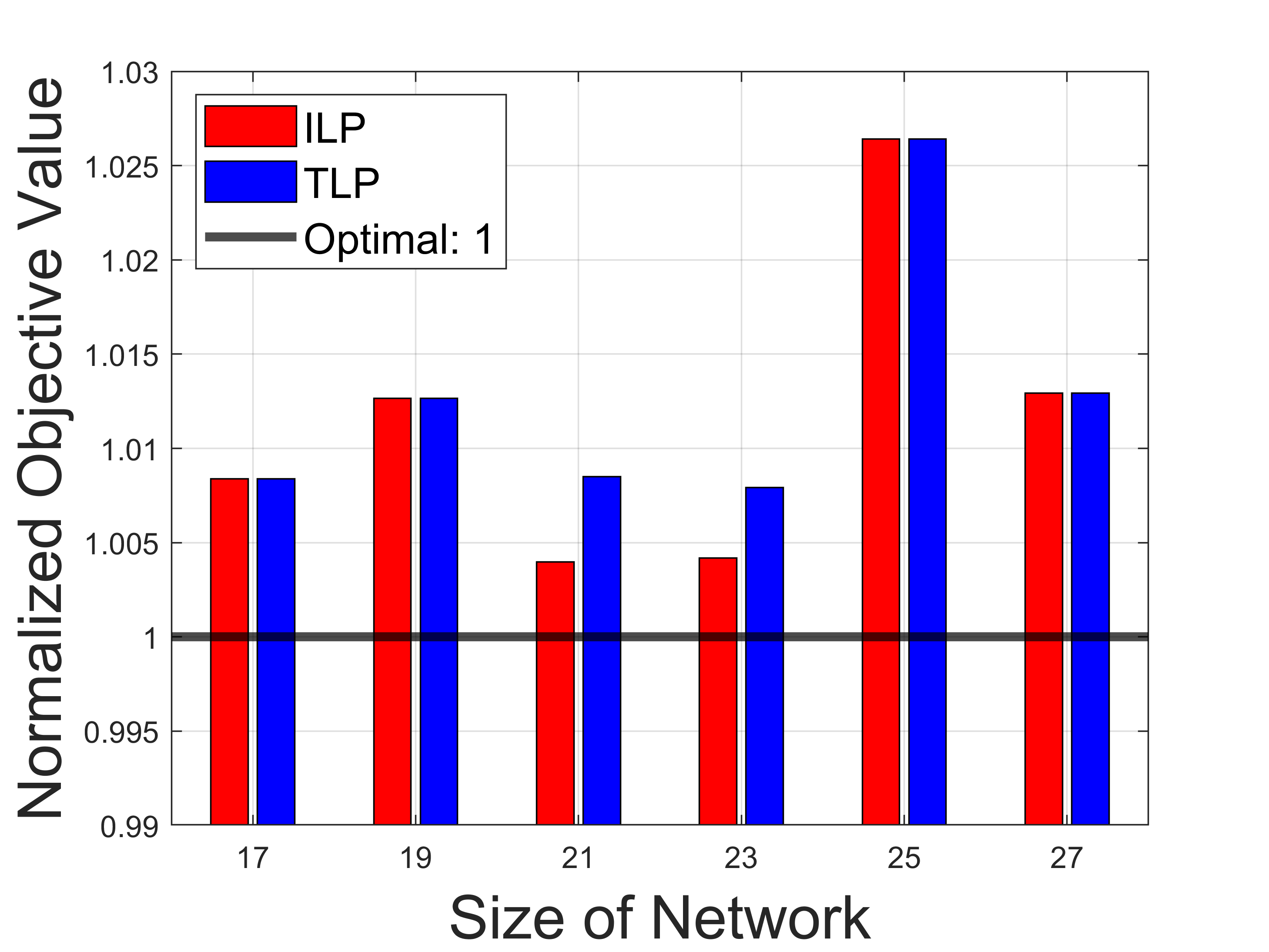}
\caption{\textcolor{black}{ Objective value of road map of Belgium}}
\label{ObjBelgium}
\end{figure}

\begin{figure}
\centering
\includegraphics[width=.65\linewidth]{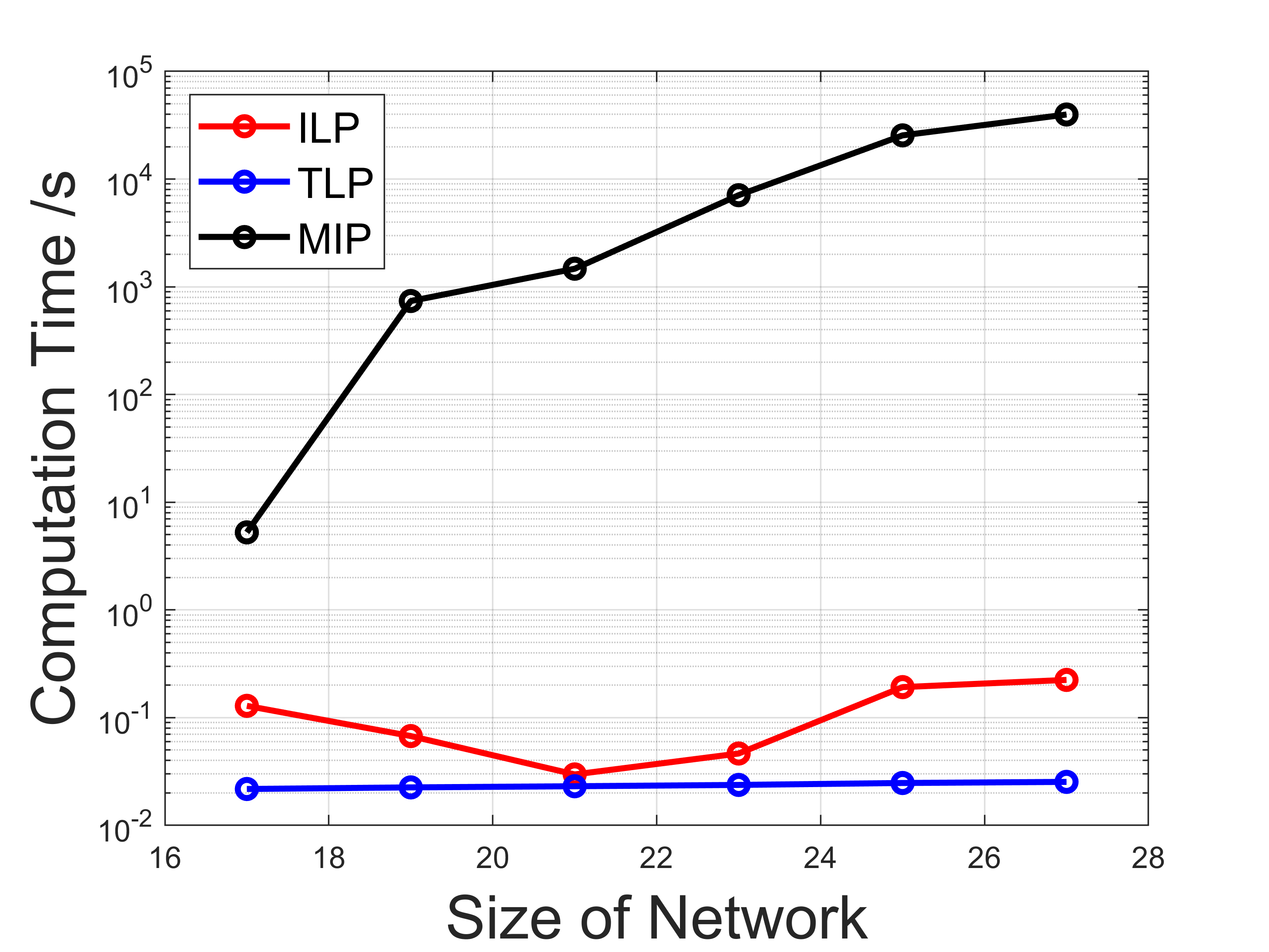}
\caption{\textcolor{black}{ Computation time of road map of Belgium}}
\label{SpeedBelgium}
\end{figure}

 \textcolor{black}{ As shown in Fig.\ref{ObjBelgium} and Fig.\ref{SpeedBelgium}, our proposed methods applied in real system of Belgium road also achieve great performance in solution quality and computation time reduction, which validates the effectiveness and superiority of proposed algorithms.
}

\section{Conclusion and Discussion}
In this paper, we investigate the joint routing and charging problem of multiple logistic EVs, which is formulated as an MIP problem. However, directly solving the MIP problem is very time consuming and even impractical for large scale problems in practice. To solve this problem efficiently, we first propose a two-stage optimization scheme, which solves the routing problem and charging problem sequentially. Meanwhile, by leveraging the techniques of exact LP relaxation and bilinear term elimination, both routing and charging problems can be solved in an LP fashion. Then, an iterative variant method of the two-stage one is proposed to further improve the quality of solution. Extensive simulations have been performed to verify the effectiveness of the proposed methods. 

\par\textcolor{black}{  The proposed algorithms can be applied in practice nicely, except for the following minor points. In stead of global optimal solution, the near-optimal solution can be obtained by proposed methods within polynomial time. Meanwhile, large-scale dispatch of EVs may have impact on the road condition, which has not been considered in our mathematical model. In the general case, however, the impact of EVs dispatch can be ignored since the logistic EV is only a small part of the transportation system. In addition, the energy consumption relation with travel condition has been simplified in this study, which however can be easily recovered to more realistic cases when necessary.
}
\par\textcolor{black}{
Meanwhile, in some real world applications, the stochastic and time varying nature of travel condition and dynamic electricity price may also be considered. Since these parameters are static in this paper, our algorithms need some nontrivial modifications before application. A good way is to put our algorithms into the framework of rolling horizon \cite{zhu2013energy}, which is designed to solve real-time optimization problems.  In rolling horizon framework, system parameters should be updated every small time step to capture the real-time change, then we re-optimize the problem with new parameters. Clearly, the rolling horizon framework demands a very fast computation speed, which is exactly the advantage of our algorithms. However, the location of EVs will change in the re-optimization stage since they are already sent to requests, which poses a new challenge to the implementation of our algorithms, i.e., multiple start depots. By introducing a virtual depot connected with the multiple start depots via zero travel distance, the multiple start depots routing problem can be transformed as the single depot case. Then, our algorithms can apply.
}



%


\appendices




\ifCLASSOPTIONcaptionsoff
  \newpage
\fi



%


\bibliographystyle{IEEEtran}  
\bibliography{references}

%
\begin{IEEEbiography}[{\includegraphics[width=1in,height=1.25in,clip,keepaspectratio]{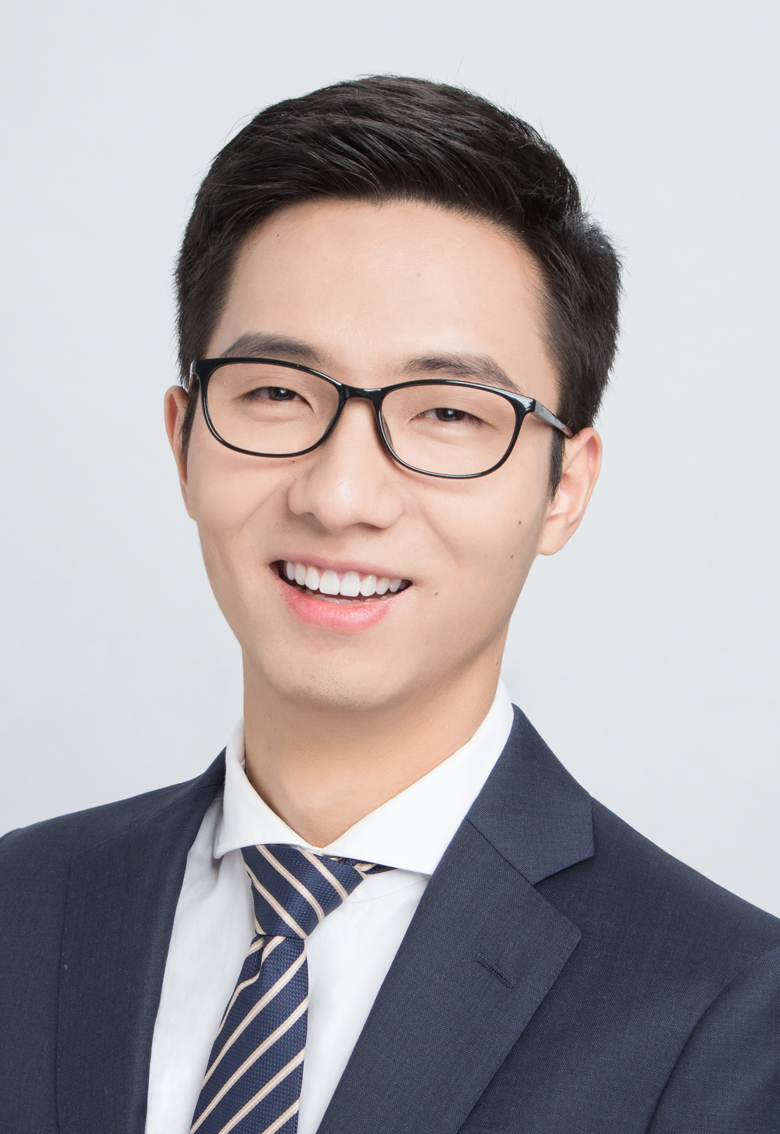}}]{Canqi Yao}
received the B.Eng. degree in electrical engineering from Changsha University of Science and Technology, Changsha, Hunan, China, in 2018. He is currently pursuing the Ph.D. degree in mechanical engineering from the Southern University of Science and Technology (SUSTech), Shenzhen, Guangdong, China. His current research interests include smart grid, electric transportation systems, and optimization theory.\end{IEEEbiography}
\begin{IEEEbiography}[{\includegraphics[width=1in,height=1.25in,clip,keepaspectratio]{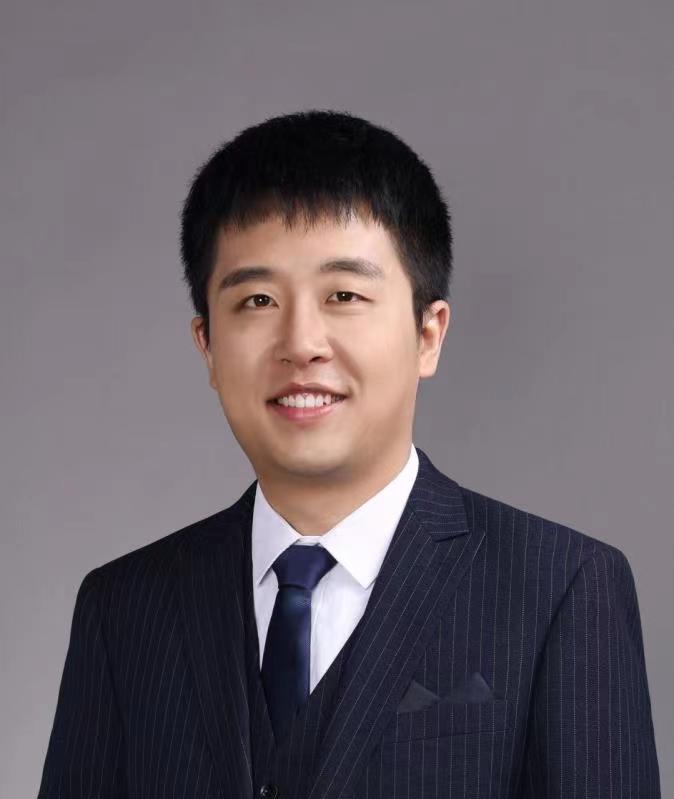}}]{Shibo Chen}
(M’17) received the B.Eng. degree in electronic engineering from University of Science and Technology of China (USTC), Hefei, China, in 2011 and the Ph.D. degree in electronic and computer engineering from the Hong Kong University of Science and Technology (HKUST), Kowloon, Hong Kong, in 2017. He was a Postdoctoral Fellow with HKUST before joining the Department of Mechanical and Energy Engineering , Southern University of Science and Technology (SUSTech), Shenzhen, China in 2019 as a Research Assistant Professor. His current research interests include smart grid, optimization theory and game theory.\end{IEEEbiography}
\begin{IEEEbiography}[{\includegraphics[width=1in,height=1.25in,clip,keepaspectratio]{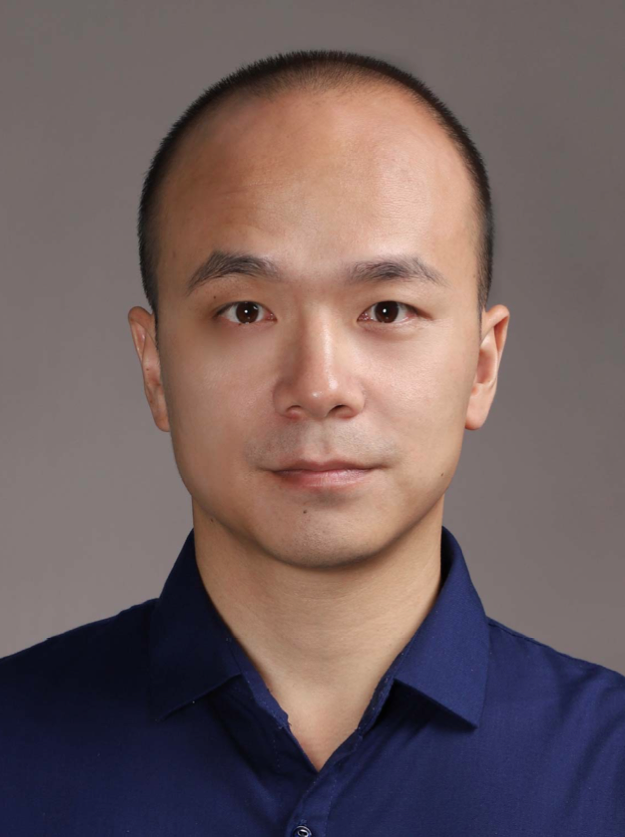}}]{Zaiyue Yang}
 (M'10) received the B.S. and M.S. degrees from the Department of Automation, University of Science and Technology of China, Hefei, China, in 2001 and 2004, respectively, and the Ph.D. degree from the Department of Mechanical Engineering, University of Hong Kong, in 2008. He was a Postdoctoral Fellow and Research Associate with the Department of Applied Mathematics, Hong Kong Polytechnic University, before joining the College of Control Science and Engineering, Zhejiang University, Hangzhou, China, in 2010. Then, he joined the Department of Mechanical and Energy Engineering, Southern University of Science and Technology, Shenzhen, China, in 2017. He is currently a Professor there. His current research interests include smart grid, signal processing and control theory. Prof. Yang is an associate editor for the IEEE Transactions on Industrial Informatics.\end{IEEEbiography}





\end{document}